\documentclass{article}%
\usepackage{amsmath}
\usepackage{amsfonts}
\usepackage{amssymb}
\usepackage{graphicx}
\usepackage{rotating}
\usepackage[ruled,vlined,linesnumbered,noresetcount]{algorithm2e}
\usepackage{hyperref}
\usepackage{colortbl}
\usepackage[round]{natbib}
\providecommand{\U}[1]{\protect\rule{.1in}{.1in}}
\newtheorem{theorem}{Theorem}

\newtheorem{lemma}[theorem]{Lemma}

\newenvironment{proof}[1][Proof]{\noindent\textbf{#1.} }{\ \rule{0.5em}{0.5em}}
\begin{document}

\title{Optimizing COVID-19 vaccine distribution adding spatio-temporal criteria}
\author{Jos\'{e} Luis Sainz-Pardo$^{1}$, Jos\'{e} Valero$^{2}$\\{\small Centro de Investigaci\'{o}n Operativa, Universidad Miguel
Hern\'{a}ndez de Elche,}\\{\small Avda. Universidad s/n, 03202, Elche (Alicante), Spain}\\$^{1}${\small jlsainz@umh.es}, $^{2}${\small jvalero@umh.es}\\}
\date{}
\maketitle

\begin{abstract}
Massive vaccination against pandemics such as Coronavirus SARS-CoV-2 presents
several complexities. The criteria to assess public health policies are
fundamental to distribute vaccines in an effective way in order to avoid as
many infections and deaths as possible. Usually these policies are focused on
determining socio-demographic groups of people and establishing a vaccination
order among these groups.

This work focuses on optimizing the way of distributing vaccines among the
different populations of a region for a period of time once established the
priority socio-demographic groups. For this aim we use a SEIR\ model which
takes into account vaccination. Also, for this model we prove theoretical results concerning the convergence
of solutions on the long-term and the stability of fixed points and analyze
the impact of an hypothetical vaccination during the COVID-19 pandemics in Spain.

After that, we introduce a heuristic approach in order to minimize the COVID-19 spreading by
planning effective vaccine distributions among the populations of a region over a period of time. As an
application, the impact of distributing vaccines in the Valencian Community (Spain) according to this method
is computed in terms of the number of saved infected individuals.

\end{abstract}

\bigskip


\bigskip


\section{Introduction}

The development of effective vaccines is a solution to the public health
crisis originated by an epidemic. This has been indeed the case of the
coronavirus disease COVID-19. The impact of the pandemic around the world
motivated an unprecedented massive and urgent vaccination. For this, it was
important to set up vaccination strategies that consider different aspects,
rules and policies. At this respect, the President of the European Commission,
Ursula von der Leyen, said: "With our Vaccination Strategy, we are helping EU
countries prepare their vaccination campaigns: who should be vaccinated first,
how to have a fair distribution and how to protect the most vulnerable. If we
want our vaccination to be successful, we need to prepare now" \cite{ursula}.
The European vaccination strategy was based on establishing priority groups
for accessing to the COVID-19 vaccines. According to \cite{eustrategy}, these
groups had been setted up by two criteria: to protect the most vulnerable
groups and individuals, and to slow down and eventually stop the spread of the
disease. The priority groups finally established in European Community were:

\begin{itemize}
\item Health care and long-term care facility workers.

\item People above 60 years of age.

\item Vulnerable population due to chronic diseases, co-morbidities and other
underlying conditions.

\item Essential workers outside the health sector.

\item Communities unable to physically distance.

\item Workers unable to physically distance.

\item Vulnerable socio-economic groups and other groups at higher risk.
\end{itemize}

The main purposes of this work are:

\begin{itemize}
\item to analyze from a theoretical point of view a SEIR\ epidemic model
taking into account vaccination;

\item to simulate the hypothetical impact of vaccination during the first wave
of the COVID-19 spreading in Spain;

\item to add the spatio-temporal component to these policies in order to
analyse the influence of different spatio-temporal vaccine distributions;

\item to develop a method to provide optimal spatio-temporal vaccine
distributions simultaneously respecting the priority groups that policies had established.
\end{itemize}

In \cite{SainzValero} we studied a modified SEIR model in which we
distinguished between detected and undetected infected people, as it is known
that at least during the first wave of the COVID-19 spreading there was a high
number of undetected infected people, and they had a huge impact on the
evolution of the pandemic. It is assumed for simplicity that infected people
who are detected are posed in quarantine and then they cannot infect other
people any more, so only undetected infected individuals are able to spread
the epidemic. It is of course more realistic to consider that the rate of
infection of those who are detected is lower, as given in \cite{LiPei}. In
\cite{SainzValero} we analyzed the impact of random testing on the spreading
of the epidemic. Now, we intend to analyze the impact of vaccination. There
are lots of models in the literature in which vaccination has been taken into
account (see e.g. \cite{BaiChenweiXu}, \cite{Brauer}, \cite{CMD},
\cite{KribsVelasco}, \cite{ZhaoJiangRegan}).

In Section \ref{Model}, we introduce our SEIR model with vaccination on the
base of the modified SEIR model given in \cite{SainzValero} considering the
situation when vaccination consists of two doses and that there is a given
percentage of immunity achieved after each dose is shotted.

In Section \ref{Qualitative}, we analyse the asymptotic behaviour of the
solutions of our model in the particular case where the coefficients are
constant. The system of equations is nevertheless non-autonomous due to the
fact that the number of dosis per day is time-dependent. However, if we assume
that the number of people which are vaccinated per day goes to zero as time
goes to infinity, then the system is asymptotically autonomous and the limit
system is just the SEIR model without vaccination. We prove that all the
solutions converges as time goes to infinity to a fixed point of the limit system.

In Section \ref{Qualitative2}, we modify the model by taking into account
births and natural deaths of the population. This is a more realistic model
when we consider a large period of time. In this case, there exist two
possible equilibria (the disease-free one and the endemic one). We analyze the
local and global stabilty of these fixed points in term of the so-called
reproductive number. This kind of analysis has been done in different
SEIR\ and SIR models (see e.g. \cite{Brauer}, \cite{CMD}, \cite{LiMuldowney},
\cite{LiWang}, \cite{Wangari}, \cite{WeiXue}, \cite{Xue}, \cite{XuaMabWanga}
among many others).

In Sections \ref{Spain} and \ref{Valencian}, we perform numerical simulations
which show the effect of an hypothetical vacunation during the first wave of
the COVID-19 pandemic in Spain and in the Spanish region called the Valencian
Community in 2020. For this aim, the parameters of the model are estimated by
means of a genetic algorithm called Differential Evolution. This technique was
firstly exposed in \cite{dev1} as an evolutionary method for optimizing
nonlinear functions. It has also been employed for estimating parameters of
infectious diseases models as SIR, SIS, SEIR, SEIS and others in \cite{dev2}.
Several versions of the technique can be found in the literature. We use the
one given in \cite{dev3}.

In Section \ref{distribution} a heuristic approach based on the proposed
modified SEIR\ model is introduced in order to optimize vaccine distributions
adding spatio-temporal criteria to the priority groups policies.

The aim of Section \ref{experience} is to measure the effectiveness of the
distribution method. It is reported the effectiveness of the proposed
heuristic approach by an extensive computational experience on the Spanish
region called Valencian Community. At this regard, the distribution for
several total number of vaccines during the period from the 1st of Juny to the
31th of December of 2020 is simulated. Our approach select individuals to be
vaccinated within each priority group according to spatio-temporal criteria.
The obtained results are compared with the random selection of individuals
within each priority group.\newline

\bigskip

\section{The model\label{Model}}

The classical SEIR\ model is the following%
\begin{equation}
\left\{
\begin{array}
[c]{c}%
\dfrac{dS}{dt}=-\dfrac{\beta}{N}SI,\\
\dfrac{dE}{dt}=\dfrac{\beta}{N}SI-\sigma E,\\
\dfrac{dI}{dt}=\sigma E-\gamma I,\\
\dfrac{dR}{dt}=\gamma I,
\end{array}
\right.  \label{SEIR}%
\end{equation}
where $N$ is the size of the population, $S\left(  t\right)  $ is the number
of the susceptible individuals to the disease, $E(t)$ is the number of exposed
people assuming that in the incubation period they do not infect anyone,
$I\left(  t\right)  $ is the number of currently infected individuals which
are able to infect other people, $R(t)$ is the number of individuals that have
been infected and then removed from the possibility of being infected again or
of spreading infection (which includes dead, recovered people and those in
quarantine or with immunity to the disease). The constant $\beta$ is the the
average number of contacts per person per time, $\gamma$ is the rate of
removal ($1/\gamma$ is the average time after which an infected individual is
removed), and $1/\sigma$ is the average time of incubation of the disease. All
these parameters are non-negative.

For our purposes we need to modify system (\ref{SEIR}) in several ways.

First, in a real situation the coefficients of the model are not constants but
functions of time. Moreover, these functions should not be continuous in
general, because in an epidemic outburst the goberments impose restrictive
measures to the population leading to a sudden change of the rate of transmission.

Second, the variable $I\left(  t\right)  $ will consist of all currently
infected individuals (not only of those able to infect)\ and a new variable
$D\left(  t\right)  $, the number of currectly infected people which are
detected, will be introduced. We will assume the ideal situation in which any
detected individual is placed in quarantine, so this person is not able to
infect anyone from that moment. Thus, the number of people with the capacity
to infect others is $I\left(  t\right)  -D\left(  t\right)  .$ Also, we do not
take into account that there could be people which are inmune, so the variable
$R\left(  t\right)  $ will contain only dead and recovered individuals but
neither those in quarantine nor immune ones.

Third, we aim to estimate dead and recovered people among the detected ones
separately, so $R\left(  t\right)  $ is splitted into three variables:

\begin{itemize}
\item $F_{1}\left(  t\right)  $:\ number of dead individuals among the
detected ones;

\item $R_{1}\left(  t\right)  $:\ number of recovered individuals among the
detected ones;

\item $L\left(  t\right)  :$ number of removed individuals among the
undetected ones.
\end{itemize}

As a first step, we consider the situation where people with symptoms and
their direct contacts are detected, but there is no a plan for massive
vaccination. The rate of detection is given by the variable $\rho\left(
t\right)  $ ($0<\rho\left(  t\right)  <1$) and, therefore, $D\left(  t\right)
=\rho(t)I(t).$

With these new variables at hand system (\ref{SEIR}) becomes:%
\begin{align}
\dfrac{dS}{dt}  &  =-\dfrac{\beta\left(  t\right)  }{N}S\left(  t\right)
(1-\rho\left(  t\right)  )I\left(  t\right)  ,\nonumber\\
\dfrac{dE}{dt}  &  =\dfrac{\beta\left(  t\right)  }{N}S\left(  t\right)
(1-\rho\left(  t\right)  )I\left(  t\right)  -\sigma E\left(  t\right)
,\nonumber\\
\dfrac{dI}{dt}  &  =\sigma E(t)-\rho\left(  t\right)  \left(  \gamma
_{1}(t)+\gamma_{2}(t)+(1-\rho\left(  t\right)  )\overline{\gamma}\left(
t\right)  \right)  I\left(  t\right)  ,\nonumber\\
\frac{dF_{1}}{dt}  &  =\gamma_{1}(t)\rho\left(  t\right)  I(t),\nonumber\\
\frac{dR_{1}}{dt}  &  =\gamma_{2}\left(  t\right)  \rho\left(  t\right)
I(t),\nonumber\\
\frac{dL}{dt}  &  =\overline{\gamma}(t)(1-\rho\left(  t\right)  )I(t),
\label{SEIR2}%
\end{align}
and the currently detected individuals are given by%
\begin{equation}
D(t)=\rho\left(  t\right)  I\left(  t\right)  . \label{Detected}%
\end{equation}
Here, $\gamma_{1}\left(  t\right)  $ is the rate of mortality of detected
people at moment $t$, $\gamma_{2}\left(  t\right)  $ stands for the rate of
recovery of detected people at moment $t$, whereas $\overline{\gamma}\left(
t\right)  $ is the rate of removal among those undetected.

In a second step, we describe the situation where a vaccination strategy is
planned in order to slow down the pandemic. For this aim we define two new
variables $\Delta_{1}\left(  t\right)  $ and $\Delta_{2}\left(  t\right)  ,$
respectively standing for the number of people that is vaccinated with the
first and the second dose at time $t$. We denote by $\pi_{1}$ and $\pi_{2}$
the percentage of immunity achieved after each dose is shotted. Then we need
to rest from the first equation the number of susceptible people that are
immunized at moment $t$. We denote by $V(t)$ the total number of people that
have been immunized due to the vaccination until the moment $t$.

We need to approximate the probability of being already immunized by natural
ways when a vaccine is applied, so people that have been immunized naturally
and after that vaccinated are not counted twice. If we choose at moment $t$ an
arbitrary person among the whole population, then the probability of being
susceptible is $S(t)/N$. If we choose an arbitrary individual among the ones
that have not been vaccinated yet, then the probability of being susceptible
is given approximately by%
\[
\frac{S(t)-\overline{V}_{1}(t)+V(t)}{N-\overline{V}_{1}(t)},
\]
where $\overline{V}_{1}(t)$ stands for the number of people that have received
at least the first dose until the moment $t$. That is, this probability is the
proportion of people that are susceptible but have not been vaccinated among
the total number of non-vaccinated individuals. When the number of vaccinated
people is low, then this probability can be approximated by $S(t)/N$. Also,
this approximation is fair if we assume that the proportion of susceptible
people among the non-vaccinated individuals is similar to the proportion of
susceptible people in the whole population. For the second dose, we choose an
arbitrary individual among people who have received the first dose but not the
second. Then the probability of being susceptible in this case is the
proportion of susceptible people among those who have received the first dose
of the vaccine but not the second. Again, we assume that this proportion is
similar to the proportion of susceptible people in the whole population, so it
is approximated by $S(t)/N$. With these assumptions the number of susceptible
people that are immunized at moment $t$ is given by%
\[
\frac{S(t)}{N}\left(  \Delta_{1}(t)\pi_{1}+\Delta_{2}(t)(\pi_{2}-\pi
_{1})\right)
\]
and the derivative of the function $V$ by%
\[
\frac{dV}{dt}=\frac{S(t)}{N}\left(  \Delta_{1}(t)\pi_{1}+\Delta_{2}(t)(\pi
_{2}-\pi_{1})\right)  .
\]
Hence, system (\ref{SEIR2}) becomes:%
\begin{align}
\dfrac{dS}{dt}  &  =-\frac{\beta\left(  t\right)  }{N}S\left(  t\right)
\left(  (1-\rho\left(  t\right)  )I\left(  t\right)  \right)  -\frac{S(t)}%
{N}\left(  \Delta_{1}(t)\pi_{1}+\Delta_{2}(t)(\pi_{2}-\pi_{1})\right)
,\nonumber\\
\frac{dE}{dt}  &  =\frac{\beta\left(  t\right)  }{N}S\left(  t\right)  \left(
(1-\rho\left(  t\right)  )I\left(  t\right)  \right)  -\sigma E\left(
t\right)  ,\nonumber\\
\frac{dI}{dt}  &  =\sigma E(t)-\rho\left(  t\right)  \left(  \gamma
_{1}(t)+\gamma_{2}(t)+(1-\rho\left(  t\right)  )\overline{\gamma}\left(
t\right)  \right)  I\left(  t\right)  ,\nonumber\\
\frac{dF_{1}}{dt}  &  =\gamma_{1}(t)\rho\left(  t\right)  I(t),\nonumber\\
\frac{dR_{1}}{dt}  &  =\gamma_{2}\left(  t\right)  \rho\left(  t\right)
I(t),\nonumber\\
\frac{dL}{dt}  &  =\overline{\gamma}(t)(1-\rho\left(  t\right)
)I(t),\label{SEIR3}\\
\frac{dV}{dt}  &  =\frac{S(t)}{N}\left(  \Delta_{1}(t)\pi_{1}+\Delta
_{2}(t)(\pi_{2}-\pi_{1})\right)  .\nonumber
\end{align}

For simplicity we could assume that the rates of death and recovery are the
same among the detected and the undetected infected people. In such a case,
$\overline{\gamma}(t)=\gamma_{1}\left(  t\right)  +\gamma_{2}\left(  t\right)
$, and $\widetilde{\gamma}_{1}\left(  t\right)  =\gamma_{1}\left(  t\right)
,\ \widetilde{\gamma}_{2}\left(  t\right)  =\gamma_{2}\left(  t\right)  $.
Thus, model (\ref{SEIR3})\ would be the following:%
\begin{align}
\dfrac{dS}{dt}  &  =-\frac{\beta\left(  t\right)  }{N}S\left(  t\right)
\left(  (1-\rho\left(  t\right)  )I\left(  t\right)  \right)  -\frac{S(t)}%
{N}\left(  \Delta_{1}(t)\pi_{1}+\Delta_{2}(t)(\pi_{2}-\pi_{1})\right)
,\nonumber\\
\frac{dE}{dt}  &  =\frac{\beta\left(  t\right)  }{N}S\left(  t\right)  \left(
(1-\rho\left(  t\right)  )I\left(  t\right)  \right)  -\sigma E\left(
t\right)  ,\nonumber\\
\frac{dI}{dt}  &  =\sigma E(t)-\left(  \gamma_{1}(t)+\gamma_{2}(t)\right)
I\left(  t\right)  ,\nonumber\\
\frac{dF_{1}}{dt}  &  =\gamma_{1}(t)\left(  \rho\left(  t\right)  I(t)\right)
,\nonumber\\
\frac{dR_{1}}{dt}  &  =\gamma_{2}\left(  t\right)  \left(  \rho\left(
t\right)  I(t)\right)  ,\nonumber\\
\frac{dL}{dt}  &  =(\gamma_{1}\left(  t\right)  +\gamma_{2}\left(  t\right)
)\left(  (1-\rho\left(  t\right)  )I(t)\right)  ,\label{SEIR4}\\
\frac{dV}{dt}  &  =\frac{S}{N}\left(  \Delta_{1}(t)\pi_{1}+\Delta_{2}%
(t)(\pi_{2}-\pi_{1})\right)  .\nonumber
\end{align}

We can generalize this model to the case where $n$ doses of the vaccine are
given. In this case, we would have:%
\begin{align}
\dfrac{dS}{dt}  &  =-\frac{\beta\left(  t\right)  }{N}S\left(  t\right)
\left(  (1-\rho\left(  t\right)  )I\left(  t\right)  \right)  -\frac{S(t)}%
{N}\left(  \Delta_{1}(t)\pi_{1}+\sum_{i=2}^{n}\Delta_{i}(t)(\pi_{i}-\pi
_{i-1})\right)  ,\nonumber\\
\frac{dE}{dt}  &  =\frac{\beta\left(  t\right)  }{N}S\left(  t\right)  \left(
(1-\rho\left(  t\right)  )I\left(  t\right)  \right)  -\sigma E\left(
t\right)  ,\nonumber\\
\frac{dI}{dt}  &  =\sigma E(t)-\left(  \gamma_{1}(t)+\gamma_{2}(t)\right)
I\left(  t\right)  ,\nonumber\\
\frac{dF_{1}}{dt}  &  =\gamma_{1}(t)\left(  \rho\left(  t\right)  I(t)\right)
,\nonumber\\
\frac{dR_{1}}{dt}  &  =\gamma_{2}\left(  t\right)  \left(  \rho\left(
t\right)  I(t)\right)  ,\nonumber\\
\frac{dL}{dt}  &  =(\gamma_{1}\left(  t\right)  +\gamma_{2}\left(  t\right)
)\left(  (1-\rho\left(  t\right)  )I(t)\right)  ,\\
\frac{dV}{dt}  &  =\frac{S}{N}\left(  \Delta_{1}(t)\pi_{1}+\sum_{i=2}%
^{n}\Delta_{i}(t)(\pi_{i}-\pi_{i-1})\right)  .\nonumber
\end{align}

\section{Qualitative analysis in the case of constant
coefficients\label{Qualitative}}

In this section, we will carry out a qualitative analysis of the behaviour of
the solutions of system (\ref{SEIR4}) when the coefficients $\beta,\rho
,\gamma_{1},\gamma_{2}$ are constant. With this assumption and putting
together the variables $F_{1},\ R_{1},\ L$ system (\ref{SEIR4}) reads as%
\begin{align}
\dfrac{dS}{dt}  &  =-\dfrac{\beta}{N}S(1-\rho)I-\frac{S}{N}\left(  \Delta
_{1}(t)\pi_{1}+\Delta_{2}(t)(\pi_{2}-\pi_{1})\right)  ,\nonumber\\
\dfrac{dE}{dt}  &  =\dfrac{\beta}{N}S(1-\rho)I-\sigma E,\nonumber\\
\dfrac{dI}{dt}  &  =\sigma E-\gamma I,\nonumber\\
\dfrac{dR}{dt}  &  =\gamma I,\label{SEIR5}\\
\frac{dV}{dt}  &  =\frac{S(t)}{N}\left(  \Delta_{1}(t)\pi_{1}+\Delta
_{2}(t)(\pi_{2}-\pi_{1})\right)  ,\nonumber
\end{align}
where $\gamma=\gamma_{1}+\gamma_{2}.$

This system is non-autonomous and we will assume that
\[
\Delta_{i}(t)\rightarrow0\text{ as }t\rightarrow+\infty\text{, }i=1,2.
\]
Hence, the system is asymptotically autonomous and the limit system is the
following:%
\begin{align}
\dfrac{dS}{dt}  &  =-\dfrac{\beta}{N}S(1-\rho)I,\nonumber\\
\dfrac{dE}{dt}  &  =\dfrac{\beta}{N}S(1-\rho)I-\sigma E,\nonumber\\
\dfrac{dI}{dt}  &  =\sigma E-\gamma I,\nonumber\\
\dfrac{dR}{dt}  &  =\gamma I,\label{SEIRLimit}\\
\frac{dV}{dt}  &  =0.\nonumber
\end{align}

This system has an infinite number of fixed points given by%
\[
\left(  S_{\infty},0,0,N-S_{\infty}-V_{\infty},V_{\infty}\right)  ,\
\]
where $0\leq S_{\infty}+V_{\infty}\leq N$, $S_{\infty},V_{\infty}\geq0$. We
will prove that every non-negative solution of system (\ref{SEIR5}) converges
as $t\rightarrow+\infty$ to a fixed point of system (\ref{SEIRLimit}).

We are interested only in considering non-negative solutions, so, first, we
will prove that is the initial datum is non-negative, then the solution
remains non-negative for every forward moment of time.

\begin{lemma}
\label{Positive}If $S_{0},\ E_{0},\ I_{0},\ R_{0},V_{0}\geq0,\ $then $S\left(
t\right)  ,\ E\left(  t\right)  ,\ I\left(  t\right)  ,\ R\left(  t\right)
,\ V\left(  t\right)  \geq0$ for all $t\geq0.$ Also:

\begin{itemize}
\item[a)] $S\left(  t\right)  >0$, for $t\geq0$, if $S_{0}>0$, whereas
$S\left(  t\right)  \equiv0$ if $S_{0}=0.$

\item[b)] If either $S_{0}>0$ or $V_{0}>0$, then $V\left(  t\right)  >0$ for
$t>0$. If $S_{0}=V_{0}=0$, then $V\left(  t\right)  \equiv0.$

\item[c)] If $E_{0}>0$, then $E(t),I(t),R(t)>0$ for $t>0.$

\item[d)] If $I_{0}>0$, then $I(t),R(t)>0$ for $t>0.$ If, moreover, $S_{0}=0$,
then $E(t)\equiv0$, whereas $S_{0}>0$ implies $E(t)>0$ for $t>0.$

\item[e)] If $R_{0}>0$, then $R(t)>0$ for $t\geq0.$
\end{itemize}
\end{lemma}

\begin{proof}
It is easy to see by a simple integration that $S\left(  t\right)  \geq0$ if
$S_{0}\geq0$ and that $S\left(  t\right)  >0$, for $t\geq0$, if $S_{0}>0$,
whereas $S\left(  t\right)  \equiv0$ if $S_{0}=0.$ By integrating the last
equation we obtain that $V(t)\geq0$ if $S_{0},V_{0}\geq0$ and the statement b)
follows as well.

If $E_{0}=I_{0}=0$, then $I(t)=E(t)=0$ for all $t\geq0$.

Assume that $E_{0}>0$. Integrating the third equation we have%
\begin{equation}
I(t)=I_{0}e^{-\gamma t}+\sigma\int_{0}^{t}e^{-\gamma\left(  t-s\right)
}E(s)ds. \label{I}%
\end{equation}
Let us prove that $E(t)>0$ for all $t\geq0$. If not, there is $t_{0}>0$ such
that $E(t_{0})=0$ and $E(t)>0$ for $t\in\lbrack0,t_{0})$. Then equality
(\ref{I}) implies that $I(t)>0$ for $t\in(0,t_{0}]$, and from the second
equation in (\ref{SEIR5}) we obtain%
\begin{equation}
E(t)=E_{0}e^{-\sigma t}+\dfrac{\beta}{N}(1-\rho)\int_{0}^{t}e^{-\sigma
(t-s)}S(s)I(s)ds>0,\text{ for }t\in\lbrack0,t_{0}], \label{E}%
\end{equation}
which is a contradiction. Using (\ref{I}) again we obtain that $I(t)>0$ for
all $t>0$. From the fourth equation in (\ref{SEIR5}) we have $R(t)>0$ for
$t>0.$

Let now $I_{0}>0$ and $E_{0}=0$. Then a similar argument as before entails
that $I(t),R(t)>0$ for $t>0$. If $S_{0}=0$, the equality in (\ref{E}) implies
that $E(t)\equiv0$. If $S_{0}>0$, then $S(t)>0$, for $t\geq0$, and the
equality in (\ref{E}) give that $E(t)>0$ for $t>0.$

Finally, if $R_{0}>0$, then
\[
R(t)=R_{0}+\gamma\int_{0}^{t}I(s)ds>0\text{ for }t\geq0.
\]

\end{proof}

\bigskip

\begin{lemma}
If $S_{0},\ E_{0},\ I_{0},\ R_{0},V_{0}\geq0,\ $then any solution of system
(\ref{SEIR5}) converges as $t\rightarrow+\infty$ to a fixed point of system
(\ref{SEIRLimit}).
\end{lemma}

\begin{proof}
We state first that $I(t)\rightarrow0$ as $t\rightarrow+\infty$. Summimg up
the first three equations and the last one we get%
\begin{equation}
\frac{dW}{dt}=-\gamma I, \label{W}%
\end{equation}
where $W(t)=S(t)+E(t)+I(t)+V(t)$. Then $W(t)$ is a non-negative, continuous,
non-increasing function, so it has a limit $W_{\infty}$ as $t\rightarrow
+\infty$. Also,
\[
W^{\prime\prime}(t)=-\gamma I^{\prime}(t)=-\gamma(\sigma E(t)-\gamma
I(t))\geq-\gamma\sigma N=-K,
\]
where $K>0$. Arguing as in \cite[Theorem 2]{SainzValero} we obtain that
$W^{\prime}(t)\rightarrow0,$ as $t\rightarrow+\infty$, which in turn implies
that $I(t)\rightarrow0,$ as $t\rightarrow+\infty.$

Further, we prove that $E(t)\rightarrow0$. Since $S\left(  t\right)  $ is
non-increasing and bounded from below by $0$ and $V(t)$ is non-decreasing and
bounded from above by $N$, it follows that $S(t)\rightarrow S_{\infty}$,
$V(t)\rightarrow V_{\infty}$. By all these convergences we have%
\[
E(t)=W(t)-I(t)-V(t)-S(t)\rightarrow W_{\infty}-V_{\infty}-S_{\infty}.
\]
We need to check that $W_{\infty}-V_{\infty}-S_{\infty}=0$. If not, then there
would be a $t_{0}>0$ such that%
\[
E(t)\geq\frac{W_{\infty}-S_{\infty}-V_{\infty}}{2}>0\text{ for all }t\geq
t_{0}.
\]
In such a case, $I(t)\rightarrow0$ implies the existence of $t_{1}\geq t_{0}$
for which%
\[
\frac{dE}{dt}\leq-\sigma\frac{W_{\infty}-S_{\infty}-V_{\infty}}{4}\text{ for
}t\geq t_{1},
\]
so%
\[
E(t)\leq E(t_{1})-\sigma\frac{W_{\infty}-S_{\infty}-V_{\infty}}{4}\left(
t-t_{1}\right)  \rightarrow-\infty,\text{ as }t\rightarrow+\infty,
\]
which is not possible as $E(t)\geq0.$
\end{proof}

\bigskip

For any $T>0$ the integrals
\[
\int_{0}^{T}\Delta_{1}(t)dt,\ \int_{0}^{T}\Delta_{2}(t)dt
\]
give the number of vaccinated people at moment $T$ with the first and second
dose, respectively. Since $\int_{0}^{T}\Delta_{1}(t)dt\leq N,\ \int_{0}%
^{T}\Delta_{2}(t)dt\leq N,$ these integrals are convergent over $\left(
0,\infty\right)  $. Denote%
\[
\overline{\Delta}_{1}=\int_{0}^{\infty}\Delta_{1}(t)dt,\ \overline{\Delta}%
_{2}=\int_{0}^{\infty}\Delta_{2}(t)dt.
\]
Then $\overline{\Delta}_{1}$ stand for the total number of people vaccinated
with at least the first dose and $\overline{\Delta}_{2}$ for the total number
of people vaccinated with the second dose.

\begin{lemma}
The limit values $\left(  S_{\infty},0,0,N-S_{\infty}-V_{\infty},V_{\infty
}\right)  $ have to satisfy the following equation%
\begin{align}
&  \frac{\beta(1-\rho)}{\gamma N}S_{\infty}-\log S_{\infty}+\frac{\beta
(1-\rho)}{\gamma N}V_{\infty}\nonumber\\
&  =\frac{1}{N}(\pi_{1}\overline{\Delta}_{1}+(\pi_{2}-\pi_{1})\overline
{\Delta}_{2})+\frac{\beta(1-\rho)}{\gamma N}(N-R_{0})-\log S_{0}%
.\label{EqPuntoFijo}%
\end{align}

\end{lemma}

\begin{proof}
Integrating in (\ref{W}) we have%
\begin{align*}
-\gamma\int_{0}^{\infty}I(t)dt  &  =S_{\infty}+V_{\infty}-\left(  S_{0}%
+E_{0}+I_{0}+V_{0}\right) \\
&  =S_{\infty}+V_{\infty}-N+R_{0}.
\end{align*}
Also, from the first equation in (\ref{SEIR5}) we obtain%
\begin{align*}
\log S_{\infty}-\log S_{0}  &  =-\frac{\beta(1-\rho)}{N}\int_{0}^{\infty
}I(t)dt-\frac{1}{N}\left(  \pi_{1}\int_{0}^{\infty}\Delta_{1}(t)dt+(\pi
_{2}-\pi_{1})\int_{0}^{\infty}\Delta_{1}(t)dt\right) \\
&  =\frac{\beta(1-\rho)}{\gamma N}(S_{\infty}+V_{\infty}-N+R_{0})-\frac{1}%
{N}(\pi_{1}\overline{\Delta}_{1}+(\pi_{2}-\pi_{1})\overline{\Delta}_{2}),
\end{align*}
so (\ref{EqPuntoFijo}) follows.
\end{proof}

\section{Qualitative analysis of a model with birth and death
rates\label{Qualitative2}}

Model (\ref{SEIR5}) is good enough when we study a epidemy in a short period
of time. However, if the period is large (it spans for a lot of years), then
we need to take into account the birth and death rates of the population. In
this case, there is no an upper bound for the total number of vaccinated people.

For simplicity, we will assume now that the number of vaccinated people per
day is constant, that is, $\Delta_{1}\left(  t\right)  \equiv\Delta
_{1},\ \Delta_{2}\left(  t\right)  \equiv\Delta_{2}$. We consider then the
following model%
\begin{align}
\dfrac{dS}{dt}  &  =\mu N-\dfrac{\beta}{N}S(1-\rho)I-\frac{S}{N}\left(
\Delta_{1}\pi_{1}+\Delta_{2}(\pi_{2}-\pi_{1})\right)  -\mu S,\nonumber\\
\dfrac{dE}{dt}  &  =\dfrac{\beta}{N}S(1-\rho)I-\sigma E-\mu E,\nonumber\\
\dfrac{dI}{dt}  &  =\sigma E-\gamma I-\mu I,\nonumber\\
\dfrac{dR}{dt}  &  =\gamma I-\mu R,\label{SEIR6}\\
\frac{dV}{dt}  &  =\frac{S}{N}\left(  \Delta_{1}\pi_{1}+\Delta_{2}(\pi_{2}%
-\pi_{1})\right)  -\mu V,\nonumber
\end{align}
and denote%
\[
p=\Delta_{1}\pi_{1}+\Delta_{2}(\pi_{2}-\pi_{1}).
\]
Since we have chosen the constant $\mu$ to be equal for the birth and death
rates, the total amount of population remains constant and equal to $N.$

This system has two equilibria. The first one is the disease-free equilibrium
given by%
\begin{align*}
P^{\ast}  &  =\left(  S^{\ast},E^{\ast},I^{\ast},R^{\ast},V^{\ast}\right) \\
&  =\left(  \frac{\mu N^{2}}{p+\mu N},0,0,0,\frac{pN}{p+\mu N}\right)  .
\end{align*}
The second one is the endemic equilibrium given by%
\[
P^{e}=\left(  S^{e},E^{e},I^{e},R^{e},V^{e}\right)  ,
\]
where%
\[
S^{e}=\frac{(\sigma+\mu)(\gamma+\mu)N}{\sigma\beta(1-\rho)},
\]%
\begin{align*}
I^{e}  &  =\frac{\mu N^{2}-S^{e}(p+\mu N)}{\beta(1-\rho)S^{e}}\\
&  =\frac{\mu N\sigma\beta(1-\rho)-(\sigma+\mu)(\gamma+\mu)(p+\mu N)}%
{(\sigma+\mu)(\gamma+\mu)\beta(1-\rho)},
\end{align*}%
\[
E^{e}=\frac{\gamma+\mu}{\sigma}I^{e},\ R^{e}=\frac{\gamma}{\mu}I^{e}%
,\ V^{e}=\frac{p}{\mu N}S^{e}.
\]

The endemic equilibrium exists when $I^{e}>0$, which happens if and only if%
\[
R_{0}:=\frac{\mu N\sigma\beta(1-\rho)}{(\sigma+\mu)(\gamma+\mu)(p+\mu N)}>1.
\]
The value $R_{0}$ is called the reproductive number.

\begin{lemma}
\label{Positive2}If $S_{0},\ E_{0},\ I_{0},\ R_{0},V_{0}\geq0,\ $then
$S\left(  t\right)  ,\ E\left(  t\right)  ,\ I\left(  t\right)  ,\ R\left(
t\right)  ,\ V\left(  t\right)  \geq0$ for all $t\geq0.$ Also:

\begin{itemize}
\item[a)] $S\left(  t\right)  ,V(t)>0$, for $t>0.$

\item[b)] If either $E_{0}>0$ or $I_{0}>0$, then $E(t),I(t),R(t)>0$ for $t>0.$

\item[d)] If $R_{0}>0$, then $R(t)>0$ for $t\geq0.$
\end{itemize}
\end{lemma}

\begin{proof}
From the first and last equations in (\ref{SEIR6}) we have%
\[
S(t)=S_{0}e^{-\int_{0}^{t}(\beta(1-\rho)I(r)/N+p/N+\mu)dr}+\mu N\int_{0}%
^{t}e^{-\int_{s}^{t}(\beta(1-\rho)I(r)/N+p/N+\mu)dr}ds>0\text{ for }t>0,
\]%
\[
V(t)=V_{0}e^{-\mu t}+\frac{p}{N}\int_{0}^{t}e^{-\mu(t-r)}S(r)dr>0\text{ for
}t>0.
\]

If $E_{0}=I_{0}=0$, then $E(t)\equiv I(t)\equiv0$.

Integrating we have%
\begin{equation}
I(t)=I_{0}e^{-(\gamma+\mu)t}+\sigma\int_{0}^{t}e^{-(\gamma+\mu)\left(
t-s\right)  }E(s)ds, \label{I2}%
\end{equation}%
\begin{equation}
E(t)=E_{0}e^{-(\sigma+\mu)t}+\dfrac{\beta}{N}(1-\rho)\int_{0}^{t}%
e^{-(\sigma+\mu)(t-s)}S(s)I(s)ds, \label{E2}%
\end{equation}%
\begin{equation}
R(t)=e^{-\mu t}R_{0}+\gamma\int_{0}^{t}e^{-\mu\left(  t-s\right)
}I(s)ds\text{ for }t\geq0. \label{R2}%
\end{equation}

If $E_{0}>0$, then arguing as in the proof of Lemma \ref{Positive} we obtain
that $E(t),I(t),R(t)>0$ for $t>0.$ If $I_{0}>0$ and $E_{0}=0$, then using
$S(t)>0$ we obtain in the same way that $I(t),E(t),R(t)>0$ for $t>0.$

Finally, from (\ref{R2}) we have that $R(t)>0$ if $R_{0}>0.$
\end{proof}

\bigskip

Let us analyze first the local stability of the equilibria in terms of
$R_{0}.$

\begin{theorem}
\label{Stable1}If $R_{0}<1$, then the disease-free equilibrium $P^{\ast}$ is
asymptotically stable.
\end{theorem}

\begin{proof}
We can take into account only the variables $S,E,I,V$, because $R=N-S-E-I-V$.
The Jacobian matrix of the system for these four variables is%
\[
J(S,E,I,V)=\left(
\begin{array}
[c]{cccc}%
-\frac{\beta(1-\rho)}{N}I-\frac{p}{N}-\mu & 0 & -\frac{\beta(1-\rho)}{N}S &
0\\
\frac{\beta(1-\rho)}{N}I & -\mu-\sigma & \frac{\beta(1-\rho)}{N}S & 0\\
0 & \sigma & -\gamma-\mu & 0\\
\frac{p}{N} & 0 & 0 & -\mu
\end{array}
\right)
\]
and at $P^{\ast}$ we have%
\[
J(S^{\ast},0,0,V^{\ast})=\left(
\begin{array}
[c]{cccc}%
-\frac{p}{N}-\mu & 0 & -\beta(1-\rho)\frac{\mu N}{p+\mu N} & 0\\
0 & -\mu-\sigma & \beta(1-\rho)\frac{\mu N}{p+\mu N} & 0\\
0 & \sigma & -\gamma-\mu & 0\\
\frac{p}{N} & 0 & 0 & -\mu
\end{array}
\right)
\]

The eigenvalues of $J$ are given by
\[
\lambda_{1}=-\frac{p}{N}-\mu,\ \lambda_{2}=-\mu,
\]%
\[
\lambda^{2}+\left(  2\mu+\sigma+\gamma\right)  \lambda+(\mu+\sigma)(\mu
+\gamma)-\frac{\sigma\beta(1-\rho)\mu N}{p+\mu N}.
\]
The eigenvalues $\lambda_{3},\lambda_{4}$ have negative real parts if%
\[
(\mu+\sigma)(\mu+\gamma)-\frac{\sigma\beta(1-\rho)\mu N}{p+\mu N}>0,
\]
that is, when $R_{0}<1$. Then $P^{\ast}$ is asymptotically stable.
\end{proof}

\bigskip

\begin{theorem}
\label{Stable2}If $R_{0}>1$, then the endemic equilibrium $P^{e}$ is
asymptotically stable.
\end{theorem}

\begin{proof}
The Jacobian matrix at $P^{e}$ is
\[
J(S^{e},E^{e},I^{e},V^{e})=\left(
\begin{array}
[c]{cccc}%
-\frac{\beta(1-\rho)}{N}I^{e}-\frac{p}{N}-\mu & 0 & -\frac{\beta(1-\rho)}%
{N}S^{e} & 0\\
\frac{\beta(1-\rho)}{N}I^{e} & -\mu-\sigma & \frac{\beta(1-\rho)}{N}S^{e} &
0\\
0 & \sigma & -\gamma-\mu & 0\\
\frac{p}{N} & 0 & 0 & -\mu
\end{array}
\right)  .
\]
The eigenvalues are $\lambda_{1}=-\mu$ and the roots of the equation%
\[
\lambda^{3}+a_{2}\lambda^{2}+a_{1}\lambda+a_{0}=0,
\]
where%
\begin{align*}
a_{0}  &  =\frac{\beta(1-\rho)}{N}(\sigma+\mu)(\gamma+\mu)I^{e},\\
a_{1}  &  =(\sigma+\gamma+2\mu)(\frac{\beta(1-\rho)}{N}I^{e}+\frac{p}{N}%
+\mu),\\
a_{2}  &  =\frac{\beta(1-\rho)}{N}I^{e}+\sigma+\gamma+3\mu+\frac{p}{N}.
\end{align*}
By the Routh-Hurwitz stability criterion, the roots have negative real part if
and only if $a_{2},a_{0}>0$ and $a_{2}a_{1}-a_{0}>0$. The assumption $R_{0}>0$
implies that $a_{0},a_{2}>0$. Also,%
\begin{align*}
&  a_{2}a_{1}-a_{0}\\
&  =\left(  \frac{\beta(1-\rho)}{N}I^{e}+\sigma+\gamma+3\mu+\frac{p}%
{N}\right)  \left(  (\sigma+\gamma+2\mu)(\frac{\beta(1-\rho)}{N}I^{e}+\frac
{p}{N}+\mu)\right) \\
&  -\frac{\beta(1-\rho)}{N}(\sigma+\mu)(\gamma+\mu)I^{e}\\
&  >0.
\end{align*}
Thus, the endemic equilibrum is asymptotically stable.
\end{proof}

\bigskip

Let us study now the global stability of the equilibria. We will consider the
invariant region%
\[
\Sigma=\{(S,E,I,R,V)\in\mathbb{R}_{+}^{5}:S+E+I+R+V=N\}.
\]

\begin{theorem}
If $R_{0}<1$, then the disease-free equilibrium $P^{\ast}$ is globally
asymptotically stable in $\Sigma.$
\end{theorem}

\begin{proof}
We will study first the variables $S,E,I$, because the first three equations
do not depend either on $R$ or $V.$

We consider first the compact positively invariant region
\[
\Omega_{\varepsilon}=\{(S,E,I)\in\mathbb{R}_{+}^{3}:S+E+I\leq N,\ S\geq
\varepsilon\},
\]
for $\varepsilon>0$. We define the function $L$ given by
\[
L(t)=S(t)-S^{\ast}\log\frac{S(t)}{S^{\ast}}+E(t)+\frac{\sigma+\mu}{\sigma
}I(t).
\]
The derivate of $L\left(  t\right)  $ is given by%
\begin{align*}
L^{\prime}(t)  &  =S^{\prime}(t)(1-\frac{S^{\ast}}{S(t)})+E^{\prime}%
(t)+\frac{\sigma+\mu}{\sigma}I^{\prime}(t)\\
&  =\left(  \mu N-\dfrac{\beta}{N}S(t)(1-\rho)I(t)-S(t)\left(  \frac{p}{N}%
+\mu\right)  \right)  \left(  1-\frac{\mu N^{2}}{p+\mu N}\frac{1}{S(t)}\right)
\\
&  +\dfrac{\beta}{N}S(t)(1-\rho)I(t)-\frac{1}{\sigma}(\sigma+\mu)(\gamma
+\mu)I(t)\\
&  =-\frac{\left(  \mu N^{2}-S(t)(p+\mu N)\right)  ^{2}}{S(t)N(p+\mu N)}%
+\frac{1}{\sigma}(\sigma+\mu)(\gamma+\mu)(R_{0}-1)I(t).
\end{align*}

Hence, if $R_{0}<1$, then $L^{\prime}(t)\leq0$ and $L^{\prime}(t)=0$ if and
only if $S=S^{\ast},\ I=0$. By Lemma \ref{Positive2} we know that if $E_{0}%
>0$, then $I(t)>0$ for $t>0$. Hence, the point $\left(  S^{\ast},0,0\right)  $
is the largest positively invariant region in $\Omega_{\varepsilon}$ in which
$L^{\prime}(t)=0$. By the Lasalle Invariance Principle \cite[Theorem
8.3.1]{Wiggins} every solution starting at $\Omega_{\varepsilon}$ converges to
the fixed point $\left(  S^{\ast},0,0\right)  $ as $t\rightarrow+\infty.$

Further, we consider the region%
\[
\Omega=\{(S,E,I)\in\mathbb{R}_{+}^{3}:S+E+I\leq N\},
\]
Lemma \ref{Positive2} implies that any solution starting at $\Omega$ satisfies
$S(t)>0$ for any $t>0$. Hence, for $t_{0}>0$ the solution belongs to
$\Omega_{\varepsilon}$ for some $\varepsilon>0$ and applying the previous
result we obtain that it converges to $\left(  S^{\ast},0,0\right)  $ as
$t\rightarrow+\infty.$

It remains to prove that $R(t)\rightarrow0$ and $V(t)\rightarrow V^{\ast
}=pN/(p+\mu N).$

Integrating the fourth equation in (\ref{SEIR6}) we have%
\[
R(t)=R_{0}e^{-\mu t}+\gamma\int_{0}^{t}e^{-\mu(t-s)}I(s)ds.
\]
For $\varepsilon>0$ there is $T_{1}(\varepsilon)>0$ such that $I(s)\leq
\varepsilon\mu/(2\gamma)$ for all $t\geq T_{1}$. Then for $t>T_{1}$ we have%
\begin{align*}
R(t)  &  \leq R_{0}e^{-\mu t}+\gamma N\int_{0}^{T_{1}}e^{-\mu(t-s)}%
ds+\frac{\varepsilon\mu}{2}\int_{T_{1}}^{t}e^{-\mu(t-s)}ds\\
&  \leq R_{0}e^{-\mu t}+\frac{\gamma N}{\mu}(e^{-\mu(t-T_{1})}-e^{-\mu
t})+\frac{\varepsilon}{2}.
\end{align*}
Thus, there is $T_{2}(\varepsilon)\geq T_{1}(\varepsilon)$ such that
$R(t)\leq\varepsilon$ for $t\geq T_{2}$. Hence, $R(t)\rightarrow0,$ as
$t\rightarrow+\infty$, and%
\[
V(t)=N-S(t)-E(t)-I(t)-R(t)\rightarrow N-S^{\ast}=V^{\ast}.
\]

\end{proof}

\bigskip

In order to study the global asymptotic stability of the endemic equilibrium,
we will follow the proof given in \cite{LiMuldowney}. As the first three
equations do not depend on the variable $V$, we will study first the behaviour
of the vector $\left(  S(t),E(t),I(t)\right)  $ in the region%
\[
\Omega=\{(S,E,I)\in\mathbb{R}_{+}^{3}:S+E+I\leq N\}.
\]

For any solution $x\left(  \text{\textperiodcentered}\right)  $ the $\omega
$-limit set is defined by%
\[
\omega(x(\text{\textperiodcentered}))=\{y:\exists t_{n}\rightarrow
+\infty\text{ such that }x(t_{n})\rightarrow y\}.
\]
Any $\omega$-limit set in the compact positively invariant set $\Omega$ is
non-empty, compact, connected and invariant \cite[Proposition \ 8.1.3]%
{Wiggins}. The elements in $\omega(x($\textperiodcentered$))$ are called
$\omega$-limit points.

\begin{lemma}
\label{Omega1}The disease-free equilibrium $\left(  \frac{\mu N^{2}}{p+\mu
N},0,0\right)  $ is the only $\omega$-limit point on the boundary of $\Omega.$
\end{lemma}

\begin{proof}
The points of the type $\left(  0,E_{\infty},0\right)  $, $E_{\infty}>0$,
cannot be $\omega$-limit points. Indeed, let $\left(  0,E_{\infty},0\right)  $
belongs to the $\omega$-limit set $\omega(x\left(  \text{\textperiodcentered
}\right)  )$ of a solution . Since $\omega(x\left(  \text{\textperiodcentered
}\right)  )$ is invariant, for any $t_{0}>0$ there exists a solution $y\left(
\text{\textperiodcentered}\right)  $ such that $y\left(  t_{0}\right)
=\left(  0,E_{\infty},0\right)  $. The initial datum $y\left(  0\right)
=\left(  S_{0},E_{0},I_{0}\right)  $ has to satisfy that either $E_{0}$ or
$I_{0}$ are positive (as, otherwise, $y(t_{0})=(S(t_{0}),0,0)$), so by Lemma
\ref{Positive2} $I(t)>0$ for $t>0$. This contradicts that $y\left(
t_{0}\right)  =\left(  0,E_{\infty},0\right)  $. By the same argument the
points $\left(  0,0,I_{\infty}\right)  $, $I_{\infty}>0$, cannot be $\omega
$-limit points as well. If $\left(  S_{\infty},E_{\infty},I_{\infty}\right)  $
satisfy $S_{\infty}+E_{\infty}+I_{\infty}=N$ and $\left(  S_{\infty}%
,E_{\infty},I_{\infty}\right)  \in\omega(x\left(  \text{\textperiodcentered
}\right)  )$, then for $t_{0}>0$ there is a solution $y\left(
\text{\textperiodcentered}\right)  $ such that $y(t_{0})=\left(  S_{\infty
},E_{\infty},I_{\infty}\right)  $. Using again Lemma \ref{Positive2} we have
that $V(t)>0,$ for $t>0$, so that $S(t_{0})+E(t_{0})+I(t_{0})=S_{\infty
}+E_{\infty}+I_{\infty}<N$. Hence, $\left(  S_{\infty},E_{\infty},I_{\infty
}\right)  $ is not an $\omega$-limit point. Finally, any orbit starting at a
point $\left(  S_{0},0,0\right)  $ moves along the $S$-axis and
$S(t)\rightarrow\frac{\mu N^{2}}{p+\mu N}.$ Hence, $\left(  \frac{\mu N^{2}%
}{p+\mu N},0,0\right)  $ is the only $\omega$-limit point on the boundary of
$\Omega.$
\end{proof}

\begin{lemma}
\label{Omega2}If $R_{0}>1$, the disease-free equilibrium $\left(  \frac{\mu
N^{2}}{p+\mu N},0,0\right)  $ cannot be the $\omega$-limit point of any orbit
starting in the interior of $\Omega$.
\end{lemma}

\begin{proof}
Define the function $L(t)=E(t)+(\sigma+\mu)I(t)/\sigma$. Its derivative is%
\begin{align*}
L^{\prime}(t)  &  =\dfrac{\beta(1-\rho)}{N}S(t)I(t)-\sigma E(t)-\mu
E(t)+\frac{\sigma+\mu}{\sigma}\left(  \sigma E-\gamma I(t)-\mu I(t)\right) \\
&  =\frac{(\sigma+\mu)(\gamma+\mu)}{\sigma}\left(  \dfrac{\sigma\beta(1-\rho
)}{N(\sigma+\mu)(\gamma+\mu)}S(t)-1\right)  I(t).
\end{align*}
Let $\left\vert S(t)-\frac{\mu N^{2}}{p+\mu N}\right\vert <\epsilon,$ where
$\epsilon>0$. Then%
\[
L^{\prime}(t)>\frac{(\sigma+\mu)(\gamma+\mu)}{\sigma}\left(  R_{0}%
-\epsilon\dfrac{\sigma\beta(1-\rho)}{N(\sigma+\mu)(\gamma+\mu)}-1\right)
I(t)>0,
\]
if $I(t)>0$ and $\epsilon$ is small enough. By Lemma \ref{Positive2} $I(t)>0$
if either $E_{0}>0$ or $I_{0}>0$, so $\left(  \frac{\mu N^{2}}{p+\mu
N},0,0\right)  $ can only be the $\omega$-limit point of an orbit starting at
an initial point with $E_{0}=I_{0}=0$, which is in the boundary of $\Omega$.
\end{proof}

\bigskip

Further, we need to study the orbital stability of periodic solutions. Let us
consider a periodic solution $q(t)=(S(t),E(t),I(t))$ with least period
$\omega>0$ and denote $\xi=\{q(t):0\leq t\leq\omega\}$. It is called orbitally
stable if for any $\varepsilon>0$ there is $\delta>0$ such that for any
solution $x(t)$ the inequality $dist(x(0),\xi)<\delta$ implies that
$dist(x(t),\xi)<\varepsilon$ for any $t\geq0$. If, moreover, there is $b>0$
such that $dist(x(t),\xi)\underset{t\rightarrow+\infty}{\rightarrow}0$
whenever $dist(x(0),\xi)<b$, then it is said to be asymptotically orbitally
stable. This orbit is asymptotically orbitally stable with asymptotic phase if
it is asymptotically orbitally stable and there is $b>0$ such that
$dist(x(0),\xi)<b$ implies the existence of $\tau$ for which $\left\vert
x(t)-q(t-\tau)\right\vert \rightarrow0$ as $t\rightarrow+\infty.$

\begin{theorem}
\label{Periodic}Any non-constant periodic solution $q\left(
\text{\textperiodcentered}\right)  $ (with least period $\omega$), if it
exists, is asymptotically orbitally stable with asymptotic phase.
\end{theorem}

\begin{proof}
It is known \cite{LiMuldowney} that a sufficient condition to get the result
is the asymptotic stability of the linear system%
\begin{equation}
\frac{dz}{dt}=J^{[2]}(q(t))z(t), \label{LinearSyst}%
\end{equation}
where $J^{[2]}$ is the second compound matrix of the Jacobian matrix
$J=(J_{ij})_{i,j\in\{1,2,3\}}$ given by%
\begin{align*}
J^{[2]}  &  =\left(
\begin{array}
[c]{ccc}%
J_{11}+J_{22} & J_{23} & -J_{13}\\
J_{32} & J_{11}+J_{33} & J_{12}\\
-J_{31} & J_{21} & J_{22}+J_{33}%
\end{array}
\right) \\
&  =\left(
\begin{array}
[c]{ccc}%
-\frac{\beta(1-\rho)}{N}I-\frac{p}{N}-\sigma-2\mu & \frac{\beta(1-\rho)}{N}S &
\frac{\beta(1-\rho)}{N}S\\
\sigma & -\frac{\beta(1-\rho)}{N}I-\frac{p}{N}-\gamma-2\mu & 0\\
0 & \frac{\beta(1-\rho)}{N}I & -2\mu-\gamma-\sigma
\end{array}
\right)  .
\end{align*}

By Lemma \ref{Omega2} the orbit $\xi$ of the periodic solution $q\left(
t\right)  =(S(t),E(t),I(t))$ remains at a positive distance from the boundary
of $\Omega$. Thus, the function%
\[
V(t)=\sup\{\left\vert X(t)\right\vert ;\frac{E(t)}{I(t)}(\left\vert
Y(t)\right\vert +\left\vert Z(t)\right\vert )\}
\]
is well defined and continuous, where $\left(  X(t),Y(t),Z(t)\right)  $ is a
solution to system (\ref{LinearSyst}). The right-hand derivative $D_{+}V(t)$
of $V(t)$ exists and a direct calculation shows that%
\begin{align*}
D_{+}\left\vert X(t)\right\vert  &  \leq-\left(  \frac{\beta(1-\rho)}%
{N}I(t)+\frac{p}{N}+\sigma+2\mu\right)  \left\vert X(t)\right\vert \\
&  +\frac{\beta(1-\rho)}{N}S(t)(\left\vert Y(t)\right\vert +\left\vert
Z(t)\right\vert ),
\end{align*}%
\[
D_{+}\left\vert Y(t)\right\vert \leq\sigma\left\vert X(t)\right\vert -\left(
\frac{\beta(1-\rho)}{N}I(t)+\frac{p}{N}+\gamma+2\mu\right)  \left\vert
Y(t)\right\vert ,
\]%
\[
D_{+}\left\vert Z(t)\right\vert \leq\frac{\beta(1-\rho)}{N}I(t)\left\vert
Y(t)\right\vert -(\gamma+\sigma+2\mu)\left\vert Z(t)\right\vert ,
\]
so%
\begin{align*}
D_{+}\frac{E(t)}{I(t)}(\left\vert Y(t)\right\vert +\left\vert Z(t)\right\vert
)  &  =\left(  \frac{E^{\prime}(t)}{E(t)}-\frac{I^{\prime}(t)}{I(t)}\right)
\frac{E(t)}{I(t)}(\left\vert Y(t)\right\vert +\left\vert Z(t)\right\vert )\\
&  +\frac{E(t)}{I(t)}D_{+}(\left\vert Y(t)\right\vert +\left\vert
Z(t)\right\vert )\\
&  \leq\sigma\frac{E(t)}{I(t)}\left\vert X(t)\right\vert +\left(
\frac{E^{\prime}(t)}{E(t)}-\frac{I^{\prime}(t)}{I(t)}-\gamma-2\mu\right)
\frac{E(t)}{I(t)}(\left\vert Y(t)\right\vert +\left\vert Z(t)\right\vert ).
\end{align*}
It follows that%
\[
D_{+}V(t)\leq\max\{g_{1}(t),g_{2}(t)\}V(t),
\]
where%
\[
g_{1}(t)=-\left(  \frac{\beta(1-\rho)}{N}I(t)+\frac{p}{N}+\sigma+2\mu\right)
+\frac{\beta(1-\rho)}{N}\frac{S(t)I(t)}{E(t)},
\]%
\[
g_{2}(t)=\sigma\frac{E(t)}{I(t)}+\frac{E^{\prime}(t)}{E(t)}-\frac{I^{\prime
}(t)}{I(t)}-\gamma-2\mu.
\]

Since
\[
\frac{\beta(1-\rho)}{N}\frac{S(t)I(t)}{E(t)}=\frac{E^{\prime}(t)}{E(t)}%
+\sigma+\mu,
\]
we have%
\[
g_{1}(t)\leq\frac{E^{\prime}(t)}{E(t)}-\mu.
\]
Also,
\[
\sigma\frac{E(t)}{I(t)}=\frac{I^{\prime}(t)}{I(t)}+\gamma+\mu
\]
gives%
\[
g_{2}(t)\leq\frac{E^{\prime}(t)}{E(t)}-\mu.
\]
Thus,%
\[
\int_{0}^{\omega}\max\{g_{1}(t),g_{2}(t)\}dt\leq\log E(\omega)-\log
E(0)-\mu\omega=-\mu\omega<0.
\]
For any $t>0$ let $n(t)$ be the largest integer such that $n(t)\omega\leq t$.
Hence, by the Gronwall lemma,%
\[
V(t)\leq V(0)e^{-\mu\omega^{2}n(t)}e^{\int_{n(t)\omega}^{t}\max\{g_{1}%
(t),g_{2}(t)\}dt}\rightarrow0,\text{ as }t\rightarrow+\infty,
\]
so $\left(  X(t),Y(t),Z(t)\right)  \rightarrow0$. Therefore, the linear system
(\ref{LinearSyst}) is asymptotically stable.
\end{proof}

\bigskip

\begin{theorem}
Any $\omega$-limit set in the interior of $\Omega$ is either a closed orbit or
the endemic equilibrium.
\end{theorem}

\begin{proof}
Our system is competitive in the interior of $\Omega$, which means that for
some diagonal matrix $H$ with values equal either to $1$ or $-1$, the matrix
$HJ(S,E,I)H$ has non-positive off diagonal elements. Indeed, for
$H=diag(-1,1,-1)$ we have%
\begin{align*}
&  H\left(
\begin{array}
[c]{ccc}%
-\frac{\beta(1-\rho)}{N}I-\frac{p}{N}-\mu & 0 & -\frac{\beta(1-\rho)}{N}S\\
\frac{\beta(1-\rho)}{N}I & -\mu-\sigma & \frac{\beta(1-\rho)}{N}S\\
0 & \sigma & -\gamma-\mu
\end{array}
\right)  H\\
&  =\left(
\begin{array}
[c]{ccc}%
-\frac{\beta(1-\rho)}{N}I-\frac{p}{N}-\mu & 0 & -\frac{\beta(1-\rho)}{N}S\\
-\frac{\beta(1-\rho)}{N}I & -\mu-\sigma & -\frac{\beta(1-\rho)}{N}S\\
0 & -\sigma & -\gamma-\mu
\end{array}
\right)  .
\end{align*}
By a version of the Poincar\'{e}-Bendixson theorem \cite[Theorem
4.1]{LiMuldowney} if $L$ is a non-empty $\omega$-limit set in the interior of
$\Omega$ and $L$ contains no equilibria, then $L$ is a closed orbit. Then if
$L$ is an $\omega$-limit set in the interior of $\Omega$ which does contain
the endemic equilibrium (which is the unique possible equilibrium in this
region), then it is a closed orbit. On the other hand, if the endemic
equilibrium $P_{1}^{e}=(S^{e},E^{e},I^{e})$ exists and $L$ contains it, as it
is asymptotically stable by Theorem \ref{Stable2}, then any solution that gets
arbitraily close to $P^{e}$ has to converge to it. Thus, $L=\{P^{e}\}.$
\end{proof}

\begin{theorem}
\label{GlobalStable2}If $R_{0}>1$, then the endemic equilibrium $P_{1}%
^{e}=(S^{e},E^{e},I^{e})$ is globally asymptotically stable in the interior of
$\Omega.$
\end{theorem}

\begin{proof}
The proof is exactly the same as in \cite[p.163]{LiMuldowney}.
\end{proof}

\begin{theorem}
If $R_{0}>1$ and the initial condition $\left(  S_{0},E_{0},I_{0},R_{0}%
,V_{0}\right)  \in\Sigma$ is such that $\left(  S_{0},E_{0},I_{0}\right)  \in
int\ \Omega$, then the solution converges to $P^{e}$ as $t\rightarrow+\infty.$
\end{theorem}

\begin{proof}
By Theorem \ref{GlobalStable2} we know that $\left(  S(t),E(t),I(t)\right)
\rightarrow P_{1}^{e}$.

Let us prove that $R(t)\rightarrow R^{e}$. In view of (\ref{R2})\ and
$R_{0}e^{-\mu t}\rightarrow0$, it is enough to show that%
\[
\gamma\int_{0}^{t}e^{-\mu\left(  t-s\right)  }I(s)ds\rightarrow R^{e}%
=\frac{\gamma}{\mu}I^{e}.
\]
We note that%
\begin{align*}
\int_{0}^{t}\mu e^{-\mu\left(  t-s\right)  }I(s)ds  &  =\int_{0}^{t}\mu
e^{-\mu\left(  t-s\right)  }I^{e}ds+\int_{0}^{t}\mu e^{-\mu\left(  t-s\right)
}\left(  I(s)-I^{e}\right)  ds\\
&  =I^{e}(1-e^{-\mu t})+\int_{0}^{t}\mu e^{-\mu\left(  t-s\right)  }\left(
I(s)-I^{e}\right)  ds,
\end{align*}
so it is enough to check that $\int_{0}^{t}\mu e^{-\mu\left(  t-s\right)
}\left(  I(s)-I^{e}\right)  ds\rightarrow0$. Since $I(t)\rightarrow I^{e}$,
for any $\varepsilon>0$ there is $T_{1}(\varepsilon)$ such that $\left\vert
I(t)-I^{e}\right\vert <\varepsilon/2$ for any $t\geq T_{1}$. Hence, there is
$T_{2}(\varepsilon)\geq T_{1}(\varepsilon)$ such that%
\begin{align*}
&  \left\vert \int_{0}^{t}\mu e^{-\mu\left(  t-s\right)  }\left(
I(s)-I^{e}\right)  ds\right\vert \\
&  \leq\left\vert \int_{0}^{T_{1}}\mu e^{-\mu\left(  t-s\right)  }\left(
I(s)-I^{e}\right)  ds\right\vert +\left\vert \int_{T_{1}}^{t}\mu
e^{-\mu\left(  t-s\right)  }\left(  I(s)-I^{e}\right)  ds\right\vert \\
&  \leq2N(e^{-\mu(t-T_{1})}-e^{-\mu t})+\frac{\varepsilon}{2}(1-e^{-\mu
(t-T_{1})})\leq\varepsilon,
\end{align*}
for $t\geq T_{2}$.

Finally,
\[
V(t)=N-S(t)-E(t)-I(t)-R(t)\rightarrow N-S^{e}-E^{e}-I^{e}-R^{e}=V^{e}.
\]

\end{proof}

\section{Application to the COVID-19 spread in Spain\label{Spain}}

In this section, we will study the effect of an hypothetical vacunation during
the first wave of the COVID-19 pandemic in Spain in 2020.

For this aim, we will use model (\ref{SEIR4}). First, we will estimate the
parameters of the model when there is no vaccination, that is, $\Delta
_{1}(t)\equiv\Delta_{2}(t)\equiv0$. Following \cite{GutierrezVarona},
\cite{Lin}, \cite{Tang} the functions $\beta\left(  t\right)  ,\ \gamma
_{1}\left(  t\right)  ,\ \gamma_{2}\left(  t\right)  $ will be piecewise
continuous functions with\ a finite number of discontinuities such that in
each interval of continuity they have the form:%
\begin{equation}
\beta\left(  t\right)  =\beta_{0}-\beta_{1}\left(  1-e^{-\alpha\left(
t-t_{0}\right)  }\right)  \label{Beta}%
\end{equation}%
\begin{equation}
\gamma_{i}(t)=\gamma_{0,i}-\gamma_{1,i}\left(  1-e^{-\alpha_{i}\left(
t-t_{0}\right)  }\right)  ,\ i=1,2. \label{Gamma}%
\end{equation}
Also, the mean value of the incubation period of the virus is about five days
(see \cite{Lauer}, \cite{RaiShukla}), so we take $\sigma=1/5$. On the other
hand, by the study of seroprevalence in Spain during the first wave of the
COVID-19 we know that at the end of May of 2020 5.2\% of the population of
Spain had been infected by the virus (which means that about 2,400,000 people
had been infected as the population in Spain is 47 millions). Since 230,000
people were detected as having caugth the infection at that moment, the
average rate of detection during the first wave of the pandemic in Spain was
approximately equal to $0.1$. Thus, we take $\rho=0.1.$

We will estimate the parameters of the model in the period from February 20,
2020 to May 17, 2020. Taking into account the different restrictions of
confinement, we split this interval into the following four subintervals: 1)
20/02-12/03; 2) 12/03-1/04; 3) 1/04-21/04; 4) 21/04-17/05.

Using the observed values of the variables $D\left(  t\right)  ,\ F_{1}\left(
t\right)  ,\ R_{1}\left(  t\right)  $, that is, the number of currently
infected, dead and recovered people which were detected, we estimate the
parameters of the functions $\beta\left(  t\right)  ,\ \gamma_{1}\left(
t\right)  $ and $\gamma_{2}\left(  t\right)  $ in each subinterval. These
values are taken from data supplied by the Spanish Health Ministry, using only
the number of infected people detected by means of a PCR\ test. The observed
value at time $t_{i}$ are denoted by $D_{i},\ F_{1i}$ and $R_{1i}$,
respectively. We optimize in each interval the square error given by%
\begin{equation}
Error=\alpha_{1}\sqrt{\sum_{i=1}^{n}\left(  D_{i}-D\left(  t_{i}\right)
\right)  ^{2}}+\alpha_{2}\sqrt{\sum_{i=1}^{n}\left(  F_{1i}-F_{1}\left(
t_{i}\right)  \right)  ^{2}}+\alpha_{3}\sqrt{\sum_{i=1}^{n}\left(
R_{1i}-R_{1}\left(  t_{i}\right)  \right)  ^{2}}, \label{Error}%
\end{equation}
where $\alpha_{1}+\alpha_{2}+\alpha_{3}=1$. We have taken $\alpha_{1}%
=\alpha_{2}=0.35,\ \alpha_{3}=0.3.$

As on the 20th of February the number of detected active infected individuals
was equal to $3$, the real number of infected subjects is given approximately
by $3/\rho=30$. At that moment there were detected neither dead nor recovered
people, and we assume that there were no removed subjects at all. Hence, the
initial value of the problem is given by:%
\[
I_{0}=30\text{, }F_{0}=0,\ R_{0}=0,\ L_{0}=0,\ S_{0}=N-I_{0}-E_{0}-F_{0}%
-R_{0}-L_{0}.
\]
The value of $E_{0}$ has to be estimated as well.

In order to make the estimation we use the genetic algorithm called the
Differential Evolution algorithm \cite{dev3}. In this method, a given array of
numbers (or \textit{genes}) containing the parameters to be estimated is
called a chromosome. The initial population \textbf{P}$_{0}$ is a set of $N$
chromosomes, which is randomly generated. Algorithm \ref{DE} describes the
procedure by which we generate a new population from the current population
\textbf{P}$_{G}$, $G=0,1,2,...$ We denote by $x_{i,G}$, $i=1.,,,N,$ the
chromosome $i$ in the population \textbf{P}$_{G}$. This sequence of
chromosomes converges to the optimal values for the parameters.

\begin{algorithm}[h]
	\setcounter{AlgoLine}{0}
	\label{DE}
	\caption{New\_Population(\textbf{P}$_G$)}
	\For{$i=1:N$}
	{
		Randomly select $r_1, r_2, r_3 \in \{1,2,...,N\}$ such that $r_1 \neq r_2 \neq r_3 \neq i$\\
		\textbf{u}$_{i,G+1} = $\textbf{x}$_{i,G} + K ($ \textbf{x}$_{r_3,G}-$ \textbf{x}$_{i,G}) +
		F ($ \textbf{x}$_{r_1,G} - $ \textbf{x}$_{r_2,G}) $ \label{desc}\\
		\textbf{if} \textbf{u}$_{i,G+1}$ is better than \textbf{x}$_{i,G}$ \textbf{then} \textbf{x}$_{i,G+1}=$\textbf{u}$_{i,G+1}$ \textbf{else}  \textbf{x}$_{i,G+1}=$ \textbf{x}$_{i,G}$\\
	}
	\textbf{return} \textbf{P}$_{G+1}$
\end{algorithm}

The values $K,F\in(0,1]$ are constants that are selected randomly. Algorithm
\ref{DE} generates one descendant for each chromosome belonging to the current
population \textbf{P}$_{G}$. The selected chromosome is replaced by its
descendant if it is better in the sense that the value of a fitness function
is lower. The algorithm is repeated until no replacement happens within a
maximal number of successive generations.

In our case, the fitness function is given by (\ref{Error}). Each chromosome
$x_{G}$ of the current population \textbf{P}$_{G}$ (we omit here the subindex
$i$) is a vector containing the parameters of the problem. In the first
period, we assume the functions $\beta\left(  t\right)  ,\ \gamma_{i}(t)$ to
be constant, so $\beta_{1}=\alpha=\gamma_{1,1}=\alpha_{1}=\gamma_{1,2}%
=\alpha_{2}=0$. Also, $E_{0}$ needs to be estimated, so $x_{G}$ contains four
parameters:%
\[
x_{G}=\left(  E_{0,G},\beta_{0,G},\gamma_{0,1,G},\gamma_{0,2,G}\right)  .
\]
In the other intervals, we estimate the nine parameters of the functions
(\ref{Beta})-(\ref{Gamma}), that is,%
\[
x_{G}=\left(  \beta_{0,G},\beta_{1,G},\alpha_{G},\gamma_{0,1,G},\gamma
_{1,1,G},\alpha_{1,G}.\gamma_{0,2,G},\gamma_{1,2,G},\alpha_{2,G}\right)  .
\]
The number of chromosomes $N$ is set to be $100$ in the first period and
\thinspace$200$ in the other ones, whereas the maximal number of iterates is
chosen to be $1000.$

The results in each interval of time are the following:

\begin{enumerate}
\item 20/02-12/03: $\beta\left(  t\right)  =\beta_{0}=1.03758,$ $\gamma
_{1}\left(  t\right)  =\gamma_{0,1}=0.0066337,\ \gamma_{2}\left(  t\right)
=\gamma_{0,2}=0.014411,$ $E_{0}=162.36331.\ $

\item 12/03-1/04: $\beta\left(  t\right)
=0.56457-0.56451(1-e^{-0.084346\left(  t-21\right)  }),$ $\gamma_{1}\left(
t\right)  =0.010016+0.0019473(1-e^{-0.11145\left(  t-21\right)  }),$
$\gamma_{2}\left(  t\right)  =0.0034428+0.082453(1-e^{-0.026258(t-21)}).$

\item 1/04-21/04: $\beta\left(  t\right)  =1.29274\ast10^{-16}%
+0.035546(1-e^{-0.84439(t-41)}),$ $\gamma_{1}\left(  t\right)
=0.0091134-0.0038616(1-e^{-0.16832(t-41)}),$ $\gamma_{2}\left(  t\right)
=0.05408-0.022434(1-e^{-0.74667(t-41)}).$

\item 21/04-17/05: $\beta\left(  t\right)  =6.33755\ast10^{-6}%
+0.031897(1-e^{-0.045468(t-61)}),$ $\gamma_{1}%
(t)=0.0040438-0.0024332(1-e^{-0.047868(t-61)}),$ $\gamma_{2}\left(  t\right)
=0.034796+0.0040778(1-e^{-0.032499(t-61)}).$
\end{enumerate}

Further, we solve system (\ref{SEIR4}) choosing a constant number of
vaccinated people per day in each period. In particular, we choose $\Delta>0$
and put $\Delta_{1}(t)$ equal to $\Delta$ through the four intervals, that is,
we assume the situation where each day during the whole first wave of the
pandemic $\Delta$ people received the first dosis of the vaccine. Since there
exists a temporal gap between the first and second dosis of the vaccine, we
put $\Delta_{2}(t)$ to be equal to $0$ in the first period, whereas
$\Delta_{2}(t)$ is equal to $\Delta$ as well in the next three periods, that
is, we suppose that from the 12-th of March each day $\Delta$ people received
the second dosis of the vaccine. Also, with resepect to the parameters
$\pi_{1},\pi_{2}$, we have chosen $\pi_{1}=0.6$ and $\pi_{2}=0.9$, that is,
60\% of people vaccinated with the first dosis become inmune, whereas after
the second one 90\% of people obtain inmunity.

In the following table we can see the predictions of the model concerning
detected deaths at the end of the whole period (i.e. on the 17-th of May) and
detected infected people at the peak of the pandemic (i.e. on the 12-th of
April) with two values of the variable $\Delta$:

\begin{center}%
\begin{tabular}
[c]{|l|l|l|}\hline
& Deaths & Infections\\\hline
Observed values & $27693$ & $98904$\\\hline
$\Delta=50000$ & $25865$ & $90723$\\\hline
$\Delta=100000$ & $24107$ & $84070$\\\hline
\end{tabular}

\end{center}

\begin{figure}[th]
\centerline{\scalebox{0.6}{\includegraphics[angle=0]{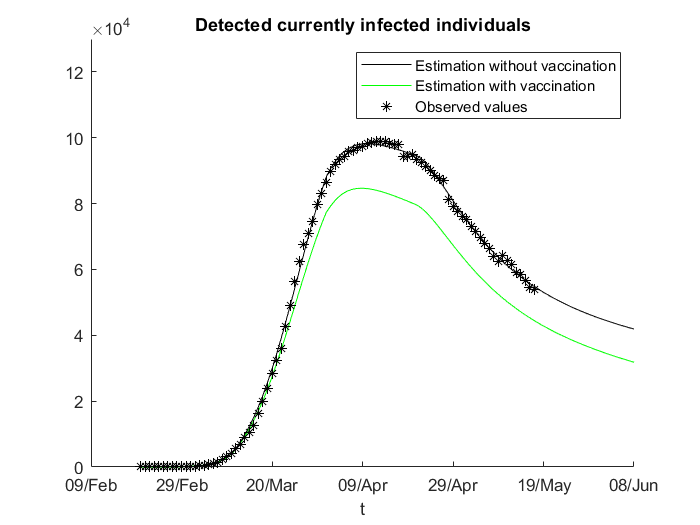}}}\caption{{\protect\footnotesize {Detected
currently infected individuals (Spain)}}}%
\label{GrafInf}%
\end{figure}

\begin{figure}[th]
\centerline{\scalebox{0.6}{\includegraphics[angle=0]{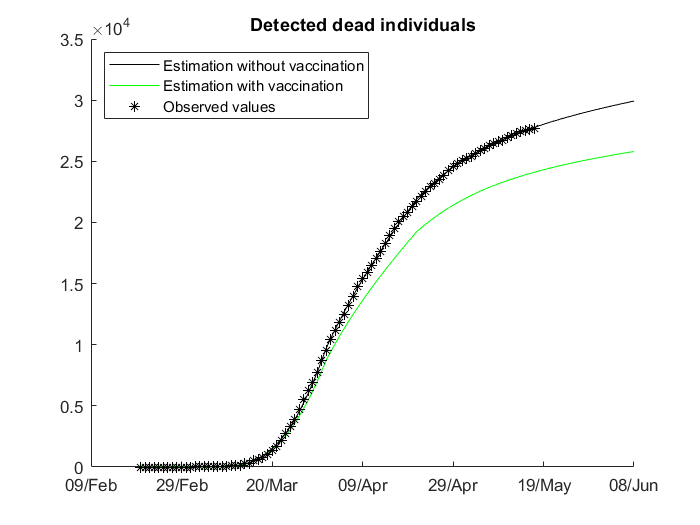}}}\caption{{\protect\footnotesize {Detected
dead individuals (Spain)}}}%
\label{GrafDead}%
\end{figure}

In Figures \ref{GrafInf}, \ref{GrafDead} we show the evolution of the number
of detected dead people and detected infected people with and without
vaccination in the case where $\Delta=100000$.

\section{Application to the COVID-19 spread in the Valencian
region\label{Valencian}}

In this section, we will study the effect of an hypothetical vacunation during
the first wave of the COVID-19 pandemic in the Spanish region called Valencian Community.

Again, we choose an average value for $\rho$ given by the study of
seroprevalence in Spain \cite{Estudio}. According to it, at the end of May of
2020 2,8\% of the population of the valencian region had been infected by the
virus (which gives about 140,000 infected people as the population is 5
millions), whereas an approximate number of 11,000 people were detected by the
COVID tests at that moment. Thus, the average rate of detection during the
first wave of the pandemic in the valencian region was approximately equal to
$0.08$.

As before, we set $\sigma=1/5.$

We estimate the parameters of the model in the period from February 25, 2020
to May 12, 2020, splitting this interval into the following five subintervals:
1) 25/02-13/03; 2) 13/03-31/03; 3) 31/03-8/04; 4) 8/04-5/05; 5) 5/05-12/05.

We proceed in the same way as in the previous section. In this case, on the
25th of February the number of detected detected active infected individuals
was $1$, so the estimate of the real number of infected subjects is
$1/\rho=13$, whereas the number of dead and recovered people, as in the case
of Spain, are equal to $0$. Hence, the initial conditions are the following%
\[
I_{0}=13\text{, }F_{0}=0,\ R_{0}=0,\ L_{0}=0,\ S_{0}=N-I_{0}-E_{0}-F_{0}%
-R_{0}-L_{0}.
\]

The estimate of the parameters in each interval is the following:

\begin{enumerate}
\item 25/02-13/03: $\beta\left(  t\right)  =\beta_{0}=0.45327,\ \gamma
_{1}\left(  t\right)  =\gamma_{0,1}=0.0047971,\ \gamma_{2}\left(  t\right)
=\gamma_{0,2}=0.0035465,\ E_{0}=122.25849.$

\item 13/03-31/03: $\beta\left(  t\right)  =2.42072-2.29381\left(
1-e^{-0.29565\left(  t-17\right)  }\right)  ,$ $\gamma_{1}\left(  t\right)
=0.016886-0.015126(1-e^{-0.048468(t-17)}),$ $\gamma_{2}\left(  t\right)
=0.0014814+0.028856(1-e^{-0.014266(t-17)}).$

\item 31/03-8/04:\ $\beta\left(  t\right)  =7.20401\ast10^{-7}-6.86704\ast
10^{-7}(1-e^{-29439.63489(t-35)}),$ $\gamma_{1}%
(t)=0.017352-0.010442(1-e^{-0.78599(t-35)}),$ $\gamma_{2}\left(  t\right)
=0.29292-0.26096(1-e^{-8.41998(t-35)}).$

\item 8/04-5/05: $\beta\left(  t\right)  =0.39963-0.38539(1-e^{-2.72216(t-43)}%
),$ $\gamma_{1}\left(  t\right)  =0.0023469+0.003184(1-e^{-1.3958(t-43)}),$
$\gamma_{2}\left(  t\right)  =0.033247+0.045749(1-e^{-0.11634(t-43)}).$

\item 5/05-12/05:\ $\beta\left(  t\right)  =\beta_{0}%
=1.834401-1.834398(1-e^{-30.03165(t-70)}),$ $\gamma_{1}%
(t)=0.0014464+0.02423(1-e^{-0.14298(t-70)}),$ $\gamma_{2}\left(  t\right)
=1.92157\ast10^{-5}+0.34632(1-e^{-1.18519(t-70)}).$
\end{enumerate}

We proceed in the same way as in the previous section by solving system
(\ref{SEIR4}) with a constant number of vaccinated people per day in each
period. Also, the distribution of the vaccines will be exactly the same,
varying only the number $\Delta$.

In the following table we can see the predictions of the model concerning
detected deaths at the end of the whole period (i.e. on the 12-th of May) and
detected infected people at the peak of the pandemic (i.e. on the 5-th of
April) with two values of the variable $\Delta$:

\begin{center}%
\begin{tabular}
[c]{|l|l|l|}\hline
& Deaths & Infections\\\hline
Observed values & $1341$ & $5767$\\\hline
$\Delta=10000$ & $1214$ & $5237$\\\hline
$\Delta=20000$ & $1102$ & $4747$\\\hline
\end{tabular}

\end{center}

\begin{figure}[th]
\centerline{\scalebox{0.6}{\includegraphics[angle=0]{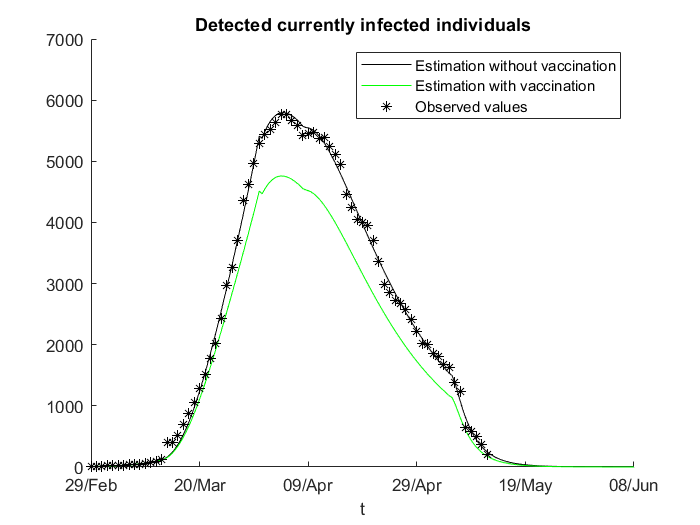}}}\caption{{\protect\footnotesize {Detected
currently infected individuals (Valencian Community)}}}%
\label{GrafInfCV}%
\end{figure}

\begin{figure}[th]
\centerline{\scalebox{0.6}{\includegraphics[angle=0]{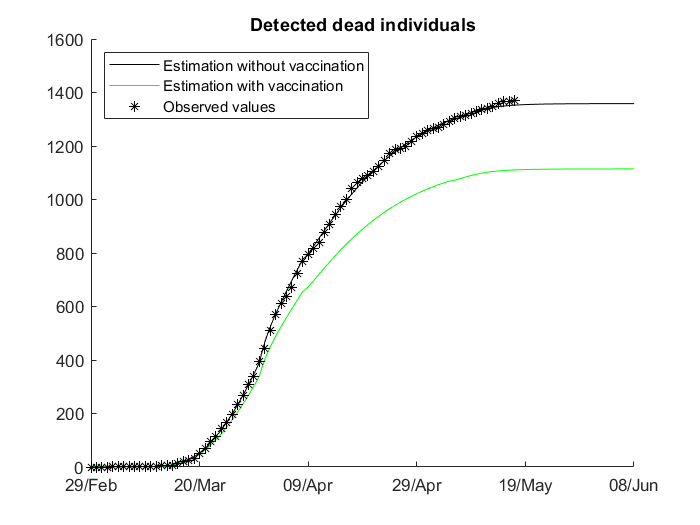}}}\caption{{\protect\footnotesize {Detected
dead individuals (Valencian Community)}}}%
\label{GrafDeadCV}%
\end{figure}

In Figures \ref{GrafInfCV}, \ref{GrafDeadCV} we can see the evolution of the
number of detected dead people and detected infected people with and without
vaccination in the case where $\Delta=20000$.

\section{Optimal distribution approach}

\label{distribution}

The aim of the distribution method is to schedule the distribution of vaccines
among the locations of a region over time within a temporal horizon but
according to the established priority groups. So, even though the space-time
component is added, no one will be vaccinated until all people belonging to
the highest priority groups have been vaccinated regardless of their location.
The addition of the spatio-temporal component is motivated by minimizing the
total number of infected cases in the region during the considered temporal
horizon. In other words, our objective is to maximize the number of saved
infected cases. We note that the saved infected cases are not only the
vaccinated people but also the people that could be infected by them in case
of contracting the disease if they were not vaccinated. Then, if we distribute
the same number of vaccines in two locations at the same instant of time, the
total number of saved infected cases due to this distribution will be
different.\newline

The distribution method is based on forecasting the number of infected people
which would be saved by vaccinating one person in one location from the first
shot to the end of the period. Then, at each instant, as many doses as
possible are assigned to the most advantageous location. However, there are
other aspects that complicate the algorithm like vaccination by priority
groups, the application of a second dose after a while, etc.\newline

\begin{algorithm}[h]
	\setcounter{AlgoLine}{0}
	\label{st1}
	\caption{gain($t$, $l$,$data$)}
	estimate parameters and state variables for $l$ from $data$ until $t-1$ \label{reup}\\
	$\bar{p} = $ forecasting of infected cases at $T$ without vaccination\\				
	$p = $ forecasting of infected cases at $T$ with $\Delta_1(t)=1$\\
	\textbf{return} $\bar{p} - p$
\end{algorithm}

Algorithm \ref{st1} forecasts infection savings due to the application of the
first shot to one person at location $l$ and instant $t$. Using the estimation
of parameters from the observed data without vaccination and the applied
vaccines until the moment of time $t-1$, we solve system (\textbf{\ref{SEIR4}%
}) and simulate new observed values of our variables by taking into account
the chosen hypothetical vaccination in this period. After that, with these new
observed data, we estimate again the parameters of the system until the moment
$t-1$ using Algorithm \textbf{\ref{DE}, }because we need to perform now the
estimation using the observed data that have been simulated taking into
account vaccination.\textbf{ }Then we solve first system (\textbf{\ref{SEIR4}%
}) until the end of the period $T$ with none vaccination, forecasting $\bar
{p}$ infected cases, and after that we solve systeml (\ref{SEIR4}) with
$\Delta_{1}(t)=1$, forecasting $p$ infected cases. We note that the day of the
application of first doses can be chosen, but not the day of the application
of the second ones. Then, the difference $\bar{p}-p$ estimates the number of
infection savings at location $l$ for the rest of the period due to the
assignation of one vaccine to location $l$ at instant $t$.\newline

\begin{algorithm}[h]
	\setcounter{AlgoLine}{0}
	\label{dist}
	\caption{optimal\_distribution($t$, $g$, $\#vaccines$, $data$)}
	\label{lim2ini}  \For {\textbf{each} $s \in rows(S)$}
	{
		\If{$S[s,2]+\delta=t$}
		{
			$d_{S[s,0],S[s,1],2}=S[s,3]$\\
			$\#vaccines=\#vaccines-S[s,3]$\\
			delete $S[s,\cdot]$\\
		}
	}\label{lim2fin}
	\While{$\#vaccines>0$}
	{				
		$l^* = $ location with maximum value of $gain(t,l,data)$ $\forall l \in L$ such that $r_{l^*, g}>0$\\
\label{dose1Ini}				$d_{l^*,g,1}=min\{r_{l^*,g}, \#vaccines\}$\\	
		$m=rows(S)+1$\\
		$S[m,0]=l^*$; $S[m,1]=g$; $S[m,2]=t$; $S[m,3]=d_{l^*,g,1}$;\\
				$\#vaccines=\#vaccines-d_{l^*,g,1}$\\				
				$r_{l^*, g}=r_{l^*, g}-d_{l^*,g,1}$\\
\label{dose1Fin}				\textbf{if} $r_{l, g}=0,$ $\forall l \in L$ \textbf{then} $g=g+1$\\					
	}
	\textbf{return} $D$
\end{algorithm}

Algorithm \ref{dist} returns the optimal matrix distribution $D=\{d_{l,g,i}\}$
which indicates the number of vaccines that will be distributed at the current
instant $t$ for each location $l$ and group $g$. The subindex $i$ is equal to
$1$ if we use our method of distribution of vaccines and equal to $2$ if we
use a random distribution. Its input parameters are the current instant $t$,
the target priority group $g$, the number of vaccines $\#vaccines$ to be
assigned at instant $t$ and the historical data of the pandemic spread.
\newline


Since $\delta$ days after the first we must assign a second dose to the people
vaccinated once, this algorithm manages a list that contains the number of
first doses which have been already applied, in which location, to which group
and at what time. That is the purpose of the list $S,$ representing its column
$S[\cdot,0]$ the location in which a first vaccination has been carried out,
whereas $S[\cdot,1]$ shows to which priority group, $S[\cdot,2]$ contains the
data of this first vaccination and $S[\cdot,3]$ contains the number of
vaccinated people. On the other hand, the matrix $\{r_{l,g}\}$ saves the
number of unvaccinated people belonging to each location $l$ and priority each
group $g$.\newline

First of all, the vaccines are assigned to the people for whom their second
vaccination deadline is $t$. We note that if the number of vaccines to be
applied every day is constant, or at least the second doses are reserved, the
distribution to this people is assured. Lines \ref{lim2ini} to \ref{lim2fin}
describe this procedure. We go through the rows of list $S$ in order to locate
which people are at the deadline date, that is, for each row $s$ such that
$S[s,2]+\delta=t$, $S[s,3]$ vaccines are assigned to the group $S[s,1]$ and
location $S[s,0]$ as second dose. Besides, the number of vaccines employed is
subtracted from the remaining amount, and since these people are fully
vaccinated, this row is removed from the list. If the number of diary
disposable vaccines is constant, according to this procedure, a period of
$\delta$ days where exclusively first doses are applied is followed by a
period of $\delta$ days where exclusively second shots are applied.\newline

In the rest of the algorithm first doses are distributed in order to achieve
the greatest effectiveness. At the current instant $t$, we find the location
$l^{\ast}$ with the highest profit after applying one first shot and such that
there are some unvaccinated people belonging to the group $g$, that is,
$r_{l^{\ast},g}>0$. Hence, we will schedule as many doses as possible to
location $l^{\ast}$ and group $g$.\newline

On the one hand, lines \ref{dose1Ini} - \ref{dose1Fin} of Algorithm \ref{dist}
show how to assign as many first doses as possible. The minimum number between
the vaccines available and the people of location $l^{\ast}$ and group $g$
pending for the first vaccination constitutes the maximum number of vaccines
to be assigned to group $g$ and location $l^{\ast}$, so $d_{l^{\ast}%
,g,1}=min\{r_{l^{\ast},g},\#vaccines\}$. Since these vaccinated people will
require a second dose, we have to add a new row to the list $S$. This is done
using the following lines of code, in which we add one row and write in its
cells the values of the location, the group, the instant and the number of
vaccines distributed. Finally, we have to subtract this number from the number
of available vaccines and from the number of unvaccinated people. If all the
people of group $g$ have already been vaccinated at all the locations then we
move on to the next target group.\newline

We also point out that, although it has only been shown how to obtain the
distribution of vaccines for instant $t,$ if obtaining a schedule in advance
for more instants is required, this can be obtained by successive projections
doing $t=t+1$ and calling Algorithm \ref{dist} again.\newline

\section{Computational experience}

\label{experience}

In this section, we will analyse the computational results of the proposed
distribution approach, comparing it with the random selection of individuals
to be vaccinated regardless of their location, but both according to the
established priority groups. First, we will explain the algorithm which is
used to simulate the number of infections, detected cases, deaths and immunity
level applying both distributions. Secondly, we will describe how the priority
groups were parametrized as well as the different instances that were
simulated. Finally, we will report the obtained results.\newline

Algorithm \ref{st3} shows the procedure which is used in order to simulate the
number of infections, detected cases, deaths and immunity level applying both
distributions. Then we estimate all the parameters and state variables for the
entire analysed period. We also set the initial instant $t_{0}$ at which the
analysed period begins and we assign the priority group to $1$.\newline

\begin{algorithm}[h]
	\setcounter{AlgoLine}{0}
	\label{st3}
	\caption{simulation}
	$\rho$ = estimate of parameters from the historical data\\				
	$t=t_0$\\
	$g_p=1$\\	
	$data_{opt}, data_{ran}$=historical data until $t_0-1$\\
	\While{$t<|T|$}
	{
		\label{ciclo1}		$O =optimal\_distribution(t,g_p,\#vaccines, data_{opt})$\\		
		$R =$ random distribution of $\#vaccines$ between priority groups\\
		$data_{opt}(t)$=simulation of data from $\rho$ at $t$ applying $O$ \\			
		$data_{prop}(t)$=simulation of data from $\rho$ at $t$ applying $R$\\			
		\label{ciclo2}		$t=t+1$\\		
	}
	\textbf{return} $\#Infections(data_{prop})$ - $\#Infections(data_{opt})$
\end{algorithm}

Since in reality there has been no vaccination during the analysed period, it
is necessary to introduce two serial of data to simulate the effect of the two
vaccine distributions for this period. Serial $data_{opt}$ will contain the
simulated infections, deaths and distributed vaccines when we distribute the
vaccines according to our approach. Serial $data_{ran}$ will be simulated by
distributing the available vaccines every day randomly selecting the
individuals of each group regardless of their location. These vaccines will be
assigned first to people with second dose deadline and the remaining to people
belonging to the priority group. Once everyone belonging to the priority group
has received their first dose then we move on to the next priority group. Note
that both serials are broken down by date and location in the same format as
the historical data. Until the instant $t_{0}-1$ both lists are equal to the
historical data.\newline

Lines \ref{ciclo1} to \ref{ciclo2} constitute a loop in which the period is
successively traversed obtaining for each instant two vaccine distributions:
the distribution using our approach $O$ and the random distribution $R$. After
that, for each instant the effect of both vaccine distributions is simulated,
obtaining new infections and deaths for each location according to each
distribution. We note that $data_{opt}$ and $data_{prop}$ are obtained
applying model (\ref{SEIR4}) with the parameters being estimated in the model
without vaccination using the historical observed data. However, each time
procedure \textit{gain} is called from $optimal\_distribution$ a new estimate
of all the parameters is carried out, as shown before in Algorithm \ref{st1}.
All the simulations are done using system (\ref{SEIR4}) with constant
coefficients.\newline

As the conditions and parameters are different according to different
circumstances as the lockdown phases, the use of masks, etc. the parameters
are defined and estimated piecewise for each week. This is similar to what is
done in Sections \ref{Spain}, \ref{Valencian} but now the intervals are chosen
in periods of seven days because in this way it can be automatically fixed,
avoiding thus the need to define the intervals depending on the government
restrictions as lockdown or curfew.\newline

Regarding the settings of the Differential Evolution technique which have been
used in this application, a population of $N=100$ individuals has been
established. The coefficients $K,F\in(0,1]$ were randomly generated at each
iteration. Finally, the estimation is stopped when a maximum number of $10000$
iterations is achieved.

We have applied our method to the pandemic spread in the Spanish region called
Valencian Community during the period from the 1st of Juny to the 31th of
December of 2020. The data on detected infections and deaths have been
collected from \cite{dades}. They are disaggregated by date and town. As it is
said, our approach adds the spatio-temporal component to the priority group
that has been established by the authorities. At this respect, in our
computational experience, we have supposed for each group the priority and
proportions which are indicated in Table \ref{priority}.\newline

\begin{table}[h]
\centering
\begin{tabular}
[c]{|c|l|r|r|}\hline
\textbf{Priority} & \textbf{Group} & \textbf{Population} & \textbf{Proportion}%
\\\hline
1 & Residents of nursing homes and health workers & 50,573 & $2.00 \%
$\\\hline
2 & Over 90 years of age & 49,307 & $0.97 \% $\\\hline
3 & People from 80 to 89 years of age & 227,224 & $4.49 \% $\\\hline
4 & People from 70 to 79 years of age & 437,862 & $8.66 \% $\\\hline
5 & People from 60 to 69 years of age & 580,728 & $11.48 \% $\\\hline
6 & People from 50 to 59 years of age & 752,334 & $14.88 \% $\\\hline
7 & People from 40 to 49 years of age & 851,588 & $16.84 \% $\\\hline
8 & People from 30 to 39 years of age & 648,759 & $12.83 \% $\\\hline
9 & People from 20 to 29 years of age & 516,126 & $10.21 \% $\\\hline
10 & People from 12 to 19 years of age & 424,531 & $8.39 \% $\\\hline
\end{tabular}
\caption{Priorities, size and proportions of each group}%
\label{priority}%
\end{table}

\begin{figure}[h]
\begin{center}
\includegraphics[width=11cm]{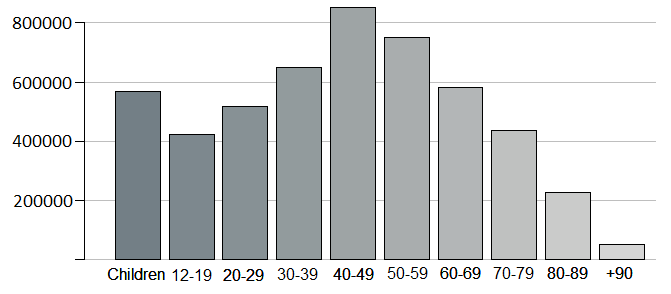}
\end{center}
\caption{Population pyramid of Valencian Community}%
\label{gPiramide}%
\end{figure}

Table \ref{priority} indicates that people over 90 years of age cannot be
vaccinated until residents of nursing homes and health workders have been
vaccinated, people from 80 to 89 years of age cannot be vaccinated until
residents of nursing homes and health workers and people over 90 have been
vaccinated, and so on. But the decision maker has to choose the people to be
vaccinated as long as they belong to the target group. Or more exactly, for
the purposes of our study, the decision maker can choose the population of
those who belong to the target group. Table \ref{priority} also shows the size
and the proportion of each group. The proportion of residents of nursing homes
and health workers has simply been assumed. However, the rest of sizes and
percentages have been extracted from \cite{ineGVA}. Children under 12 years of
age are not considered to be vaccinated, being 5,057,353 the size of the
entire population including them. Figure \ref{gPiramide} shows the population
pyramide of Valencian Community. For simplicity, we consider that these
population proportions are identical for each town, which is clearly not a
true fact. In practice, if they were known, they could be adjusted differently
for each town. Nevertheless, setting identical proportions regardless of the
place also constitutes a way for carrying out fair and equitable distributions
among the populations, which is another of the principles exposed in
\cite{ursula}. \newline

BNT162b2 is assumed to be the vaccine used for all injections. Its first dose
reaches an effectiveness of 54\% and, after the application of the second shot
21 days later, an effectiveness of 95\%. These percentages and the recommended
period between shots have been established according to \cite{efectividad}.
Regarding the recommended period, this has been strictly assumed, which means
that a second shot is applied to each vaccinated person exactly 21 days after
the first one.\newline

In addition, to get a wide computational experience, we have prepared several
examples with different number of distributed doses and percentage of people
willing to be vaccinated. These instances are:

\begin{itemize}
\item one case in which no vaccine is shot;


\item eleven cases with 10, 25, 50, 75, 100, 250, 500, 750, 1,000 thousands of doses;



\end{itemize}


We note that 1 million doses means that approximately 500 thousands people
receive their full vaccination program. The reasons for analysing up to 1,000
thousand doses are as follows:

\begin{itemize}
\item these quantities of doses are feasible to apply in the period studied;

\item a higher number than 1,000 thousand doses would imply an homogeneous
distribution, even when attempting to prioritize the populations with higher
impact, as it would require extrem mass distribution;

\item the simplification that has been proposed for the probability of being
susceptible (that is, $S(t)/N$) is valid for these quantities.
\end{itemize}


Table \ref{infectados} summarizes the obtained results relative to the
infected cases. \textit{Doses} indentifies the case, \textit{Random} and
\textit{Approach} show the estimated infected cases by using the random
distribution and our approach, \textit{R. Saving} and \textit{A. Saving}
contain the number of saved infections with both methods and, finally, Column
\textit{Advantage} shows the difference A. Saving minus R. Saving.

\begin{table}[h]
\centering
\begin{tabular}
[c]{|r|r|r|r|r|r|r|}\hline
\textbf{Doses} & \textbf{Random} & \textbf{Approach} & \textbf{R. Saving} &
\textbf{A. Saving} & \textbf{Advantage} & \\\hline
10,000 & 261,254 & 260,519 & 150 & 885 & 735 & \\\hline
25,000 & 260,856 & 259,778 & 548 & 1,626 & 1,078 & \\\hline
50,000 & 260,145 & 258,852 & 1,259 & 2,552 & 1,293 & \\\hline
75,000 & 259,831 & 258,289 & 1,573 & 3,115 & 1,542 & \\\hline
100,000 & 258,561 & 256,699 & 2,843 & 4,705 & 1,862 & \\\hline
250,000 & 255,907 & 248,866 & 5,497 & 12,538 & 7,041 & \\\hline
500,000 & 250,999 & 240,930 & 10,405 & 20,474 & 10,069 & \\\hline
750,000 & 243,438 & 230,571 & 17,966 & 30,833 & 12,867 & \\\hline
1,000,000 & 239,819 & 226,398 & 21,585 & 35,006 & 13,421 & \\\hline
\end{tabular}
\caption{Infected cases by both distribution methods}%
\label{infectados}%
\end{table}

\begin{figure}[h]
\begin{center}
\includegraphics[width=14cm]{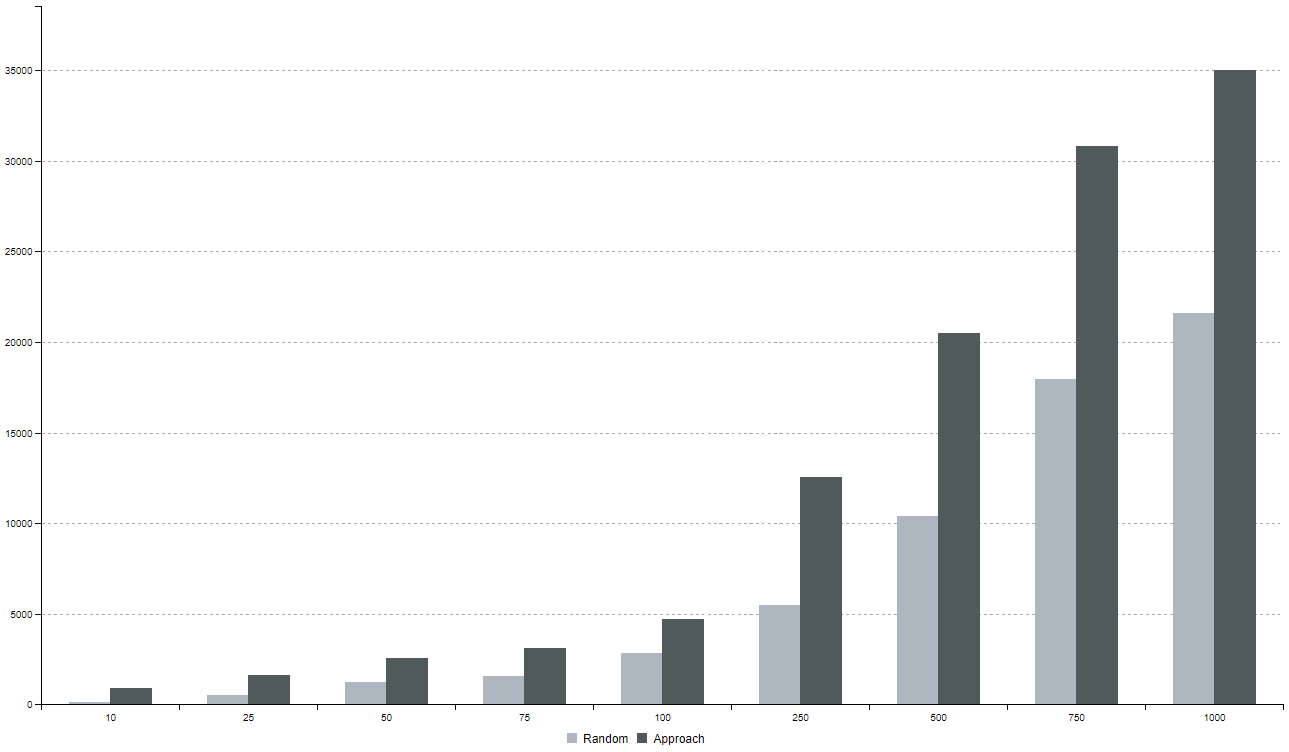}
\end{center}
\caption{Infected cases by both distribution methods}%
\label{gInfectados}%
\end{figure}

\begin{figure}[h]
\begin{center}
\includegraphics[width=12cm]{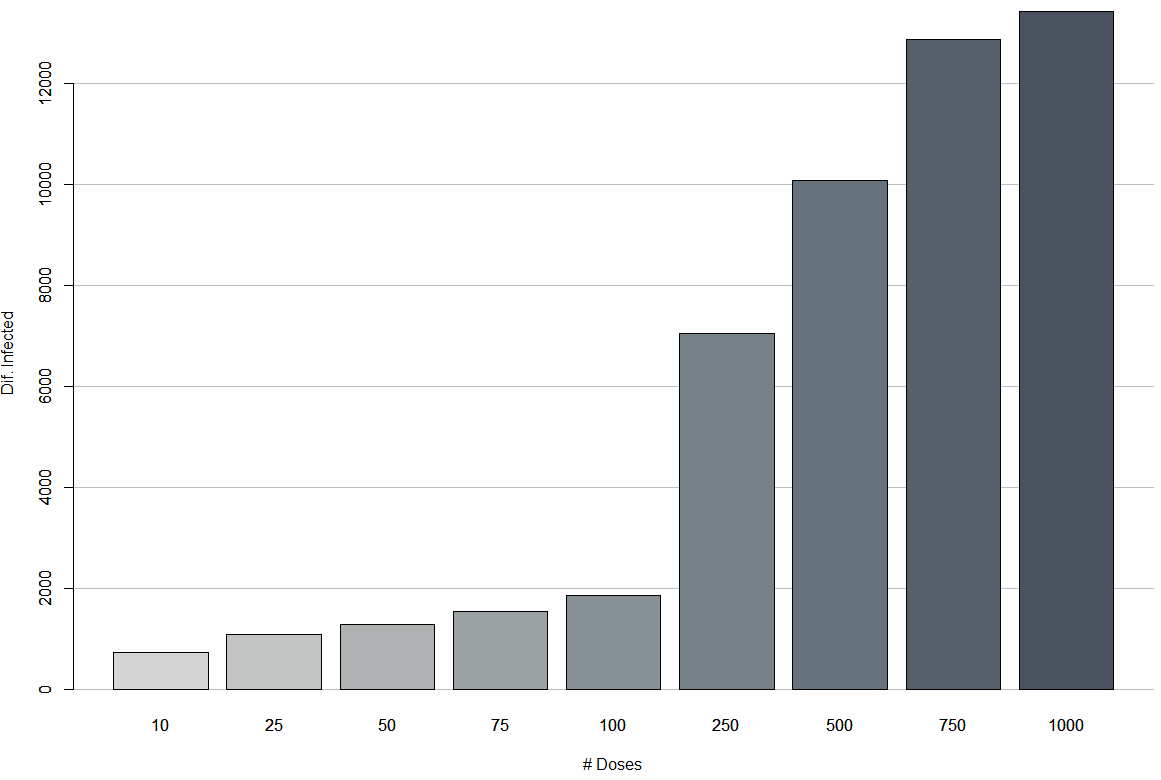}
\end{center}
\caption{Difference of infected cases by both distribution methods}%
\label{gDifInfectados}%
\end{figure}

The number of infected people was estimated to be equal to $261,404$ in the
studied period without vaccination. Hence, both saving columns are calculated
by substracting from $261,404$ the estimated values by each method. Besides,
these values are also represented in Figure \ref{gInfectados} and Column
Advantage. Finally, Figure \ref{gDifInfectados} shows the difference between
them. In both methods, it can also be seen that the more vaccines the less
infected people and, therefore, the case in which $1$ million doses are
applied is the more effective to slow down the pandemic. Besides, the proposed
method always improves the random distribution.

\begin{figure}[h]
\begin{center}
\includegraphics[width=14cm]{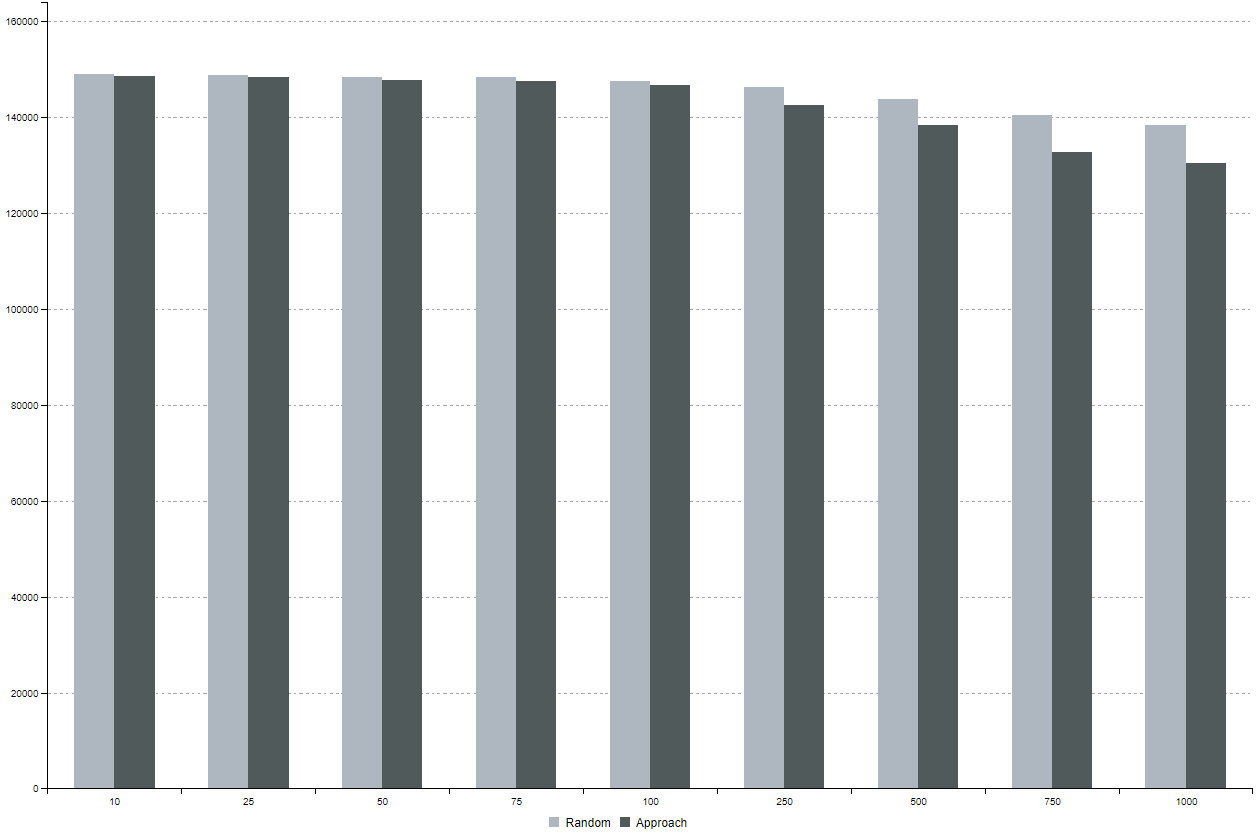}
\end{center}
\caption{Detected cases by both distribution methods}%
\label{gDetectados}%
\end{figure}

\begin{figure}[h]
\begin{center}
\includegraphics[width=12cm]{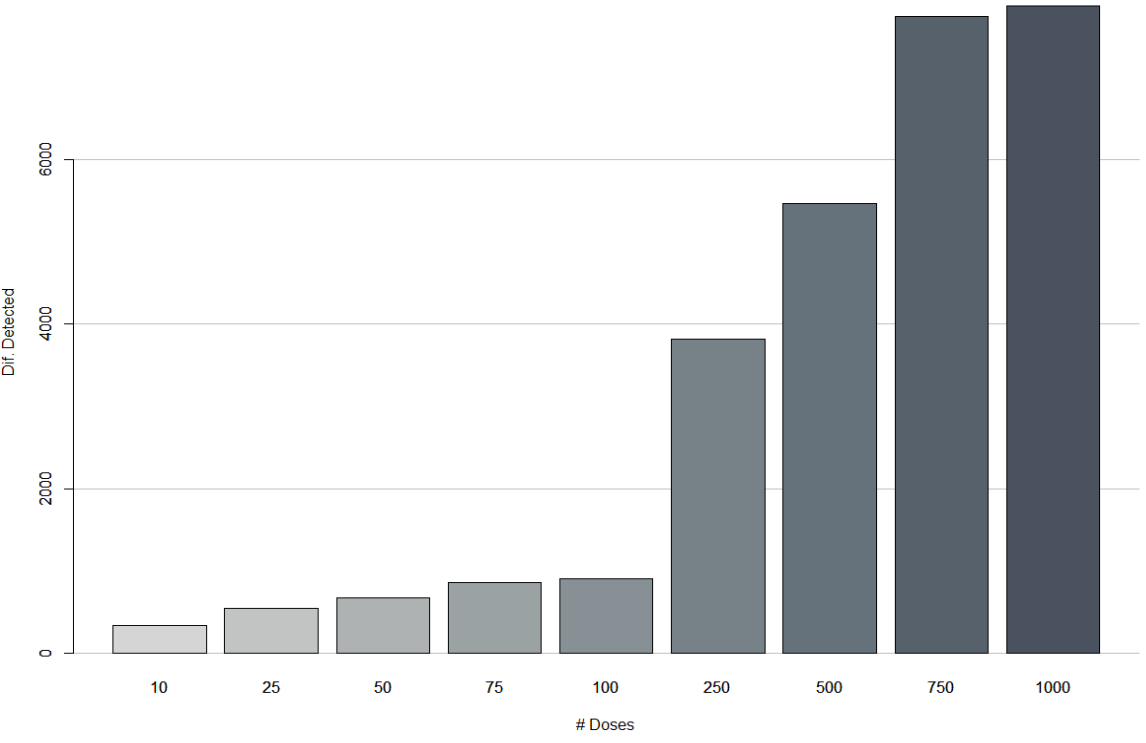}
\end{center}
\caption{Difference of detected cases by both distribution methods}%
\label{gDifDetectados}%
\end{figure}

Table \ref{detectados} and figures \ref{gDetectados} and \ref{gDifDetectados}
are relative to the number of detected cases in the simulations. They are
equivalent to the already commented table and figures of the infected cases
and the savings are calculated from the estimation of $149,754$ detected cases
for the period without vaccination. Regarding them, the results, figures and
conclusions are similar to those of the infected cases. Then, the proposed
method always improves the random distribution.\newline

\begin{table}[h]
\centering
\begin{tabular}
[c]{|r|r|r|r|r|r|}\hline
\textbf{Doses} & \textbf{Random} & \textbf{Approach} & \textbf{R. Saving} &
\textbf{A. Saving} & \textbf{Advantage}\\\hline
10,000 & 148,997 & 148,660 & 757 & 1,094 & 337\\\hline
25,000 & 148,800 & 148,262 & 954 & 1,492 & 538\\\hline
50,000 & 148,419 & 147,745 & 1,335 & 2,009 & 674\\\hline
75,000 & 148,312 & 147,452 & 1,442 & 2,302 & 860\\\hline
100,000 & 147,625 & 146,724 & 2,129 & 3,030 & 901\\\hline
250,000 & 146,295 & 142,480 & 3,459 & 7,274 & 3,815\\\hline
500,000 & 143,801 & 138,331 & 5,953 & 11,423 & 5,470\\\hline
750,000 & 140,403 & 132,664 & 9,351 & 17,090 & 7,739\\\hline
1,000,000 & 138,312 & 130,446 & 11,442 & 19,308 & 7,866\\\hline
\end{tabular}
\caption{Detected cases by both distribution methods}%
\label{detectados}%
\end{table}

Table \ref{muertes} summarizes the obtained results relative to the number of
detected deaths, and its columns and interpretation are the same as in the
previous tables. We remember that model (\ref{SEIR4}) only offers measures
about deaths that can be detected. 4,001 detected deaths are estimated in the
period without vaccination. Then, the estimated deaths are subtracted from
this value in order to get both savings. Again, as it can be guessed, the more
vaccines given, the lower the number of deaths. Besides, the proposed approach
is more advantageous in general, although the random distribution has slightly
better values in three cases with negative advantage, despite the fact that
the number of infections is always lower when our approach is used. In this
regard, it is necessary to note that the main objective of this method is to
reduce the number of infections and, as a side effect, also the number of
deaths; but the latter is not its priority.\newline

\begin{table}[h]
\centering
\begin{tabular}
[c]{|r|r|r|r|r|r|r|}\hline
\textbf{Doses} & \textbf{Random} & \textbf{Approach} & \textbf{R. Saving} &
\textbf{A. Saving} & \textbf{Advantage} & \\\hline
10,000 & 4,000 & 3,993 & 1 & 8 & 7 & \\\hline
25,000 & 3,998 & 3,971 & 3 & 30 & 27 & \\\hline
50,000 & 3,964 & 3,966 & 37 & 35 & -2 & \\\hline
75,000 & 3,978 & 3,940 & 23 & 61 & 38 & \\\hline
100,000 & 3,912 & 3,932 & 89 & 69 & -20 & \\\hline
250,000 & 3,914 & 3,824 & 87 & 177 & 90 & \\\hline
500,000 & 3,531 & 3,556 & 470 & 445 & -25 & \\\hline
750,000 & 3,779 & 3,531 & 222 & 470 & 248 & \\\hline
1,000,000 & 3,643 & 3,583 & 358 & 418 & 60 & \\\hline
\end{tabular}
\caption{Number of detected deaths by both distribution methods}%
\label{muertes}%
\end{table}

\begin{figure}[h]
\begin{center}
\includegraphics[width=14cm]{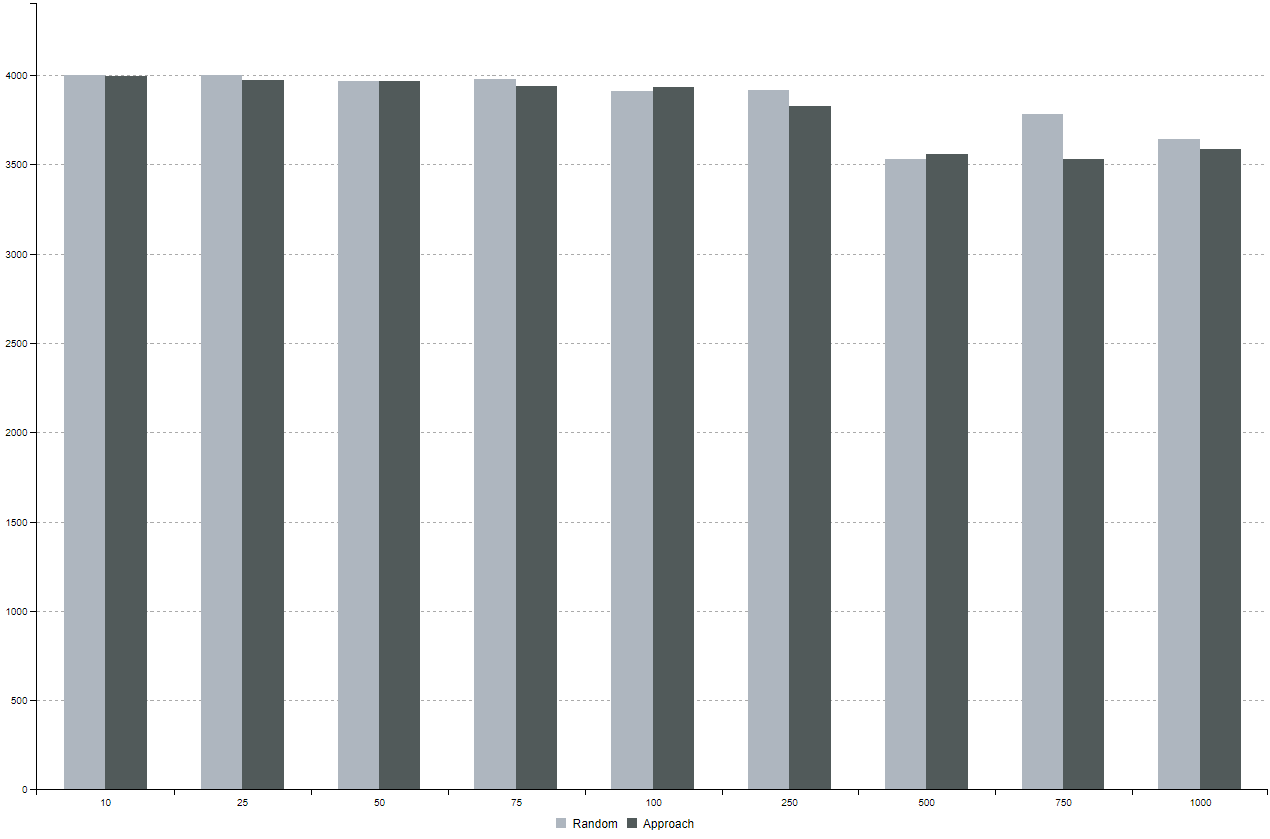}
\end{center}
\caption{Detected deaths by both distribution methods}%
\label{gMuertes}%
\end{figure}

\begin{figure}[h]
\begin{center}
\includegraphics[width=12cm]{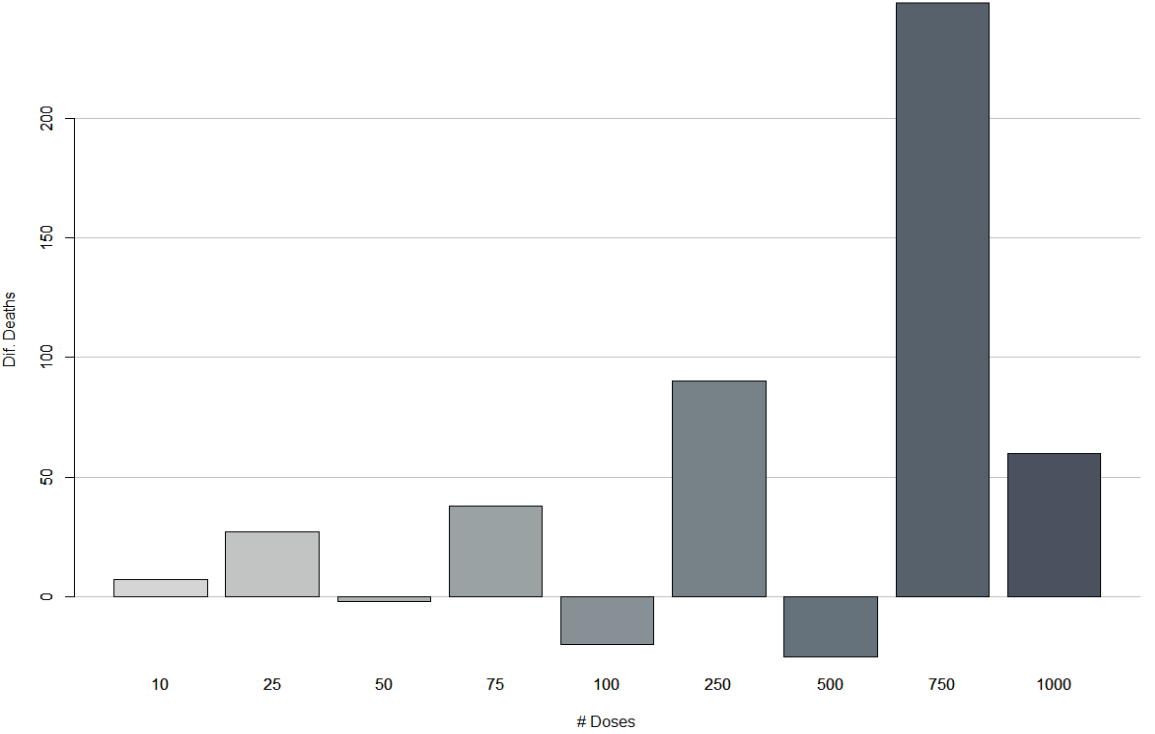}
\end{center}
\caption{Difference of detected deaths by both distribution methods}%
\label{gDifMuertes}%
\end{figure}

Figures \ref{gMuertes} and \ref{gDifMuertes} illustrate, respectively, the
number of deaths by both methods and its difference.\newline

Finally, Table \ref{inmunizados} reports the percentage of population immunity
achieved by both methods. Obviously, the more vaccination the higher
percentage of immunity. It can be observed that both methods present similar
values for the number of doses that are simulated, getting a percentage higher
than 9 \% of immunity when $1$ million vaccines are used. However, random
distribution provides a slightly higher immunisation. The explanation of this
is due to the fact that more people are infected if random distribution is
used. Therefore, it is natural immunisation that contributes to this
difference. \newline


\begin{table}[h]
\centering
\begin{tabular}
[c]{|r|r|r|r|}\hline
\textbf{Doses} & \textbf{Random} & \textbf{Approach} & \\\hline
10,000 & 8,4746 \% & 8,4631 \% & \\\hline
25,000 & 8,4841 \% & 8,4675 \% & \\\hline
50,000 & 8,5083 \% & 8,4806 \% & \\\hline
75,000 & 8,5360 \% & 8,5034 \% & \\\hline
10,0000 & 8,5489 \% & 8,5133 \% & \\\hline
25,0000 & 8,6997 \% & 8,5632 \% & \\\hline
500,000 & 9,0051 \% & 8,7607 \% & \\\hline
750,000 & 9,2016 \% & 8,9051 \% & \\\hline
100,0000 & 9,6375 \% & 9,1243 \% & \\\hline
\end{tabular}
\caption{Percentages of immunization by both distribution methods}%
\label{inmunizados}%
\end{table}


\section{Conclusions}

In this work, a SEIR model for analysing the efficiency of vaccine
distributions has been introduced. This model takes into account the presence
of both detected and non-detected infected individuals as well as the lockdown
of detected cases and, of course, vaccination. Since the values of the
parameters of the model can change abruptly due to severe governments measures
like lockdown, curfew, etc., the coefficients of the model are defined
piecewise in given intervals of time and are functions of time. This model has
been applied to the spread of the COVID\ pandemic in Spain and the Spanish
region called Valencian Community. Besides, some theoretical results
concerning convergence of solutions in the lon term and stability of
stationary points haven been proved.

On the other hand, we describe the Differential Evolution technique, which is
a genetic algorithm for the estimation of parameters, and develop a heuristic
approach for optimising the distribution of vaccines in order to minimize the
number of infected cases. This approach have been applied to the spreading of
COVID-19\ pandemic in the Valencian Community by an extensive computational
experience showing the advantages of the proposed distribution method.

\bigskip

\textbf{Acknowledgments}

This work has been partially supported by the Generalitat Valenciana (Spain), project
CIGE/2021/161. The first author has also been partially supported by the
Spanish Ministry of Science, Innovation and Universities, project
PID2021-122344NB-I00. The second author has also been partially supported by
the Spanish Ministry of Science and Innovation, projects PID2019-108654GB-I00
and PDI2021-122991NB-C21.

\bibliographystyle{apalike}
\bibliography{VacCovid}

\end{document}